\documentclass{amsart}
\usepackage{amsfonts,amssymb,mathtools,mathabx,mathrsfs, bm}
\usepackage{enumerate}

\usepackage[usenames,dvipsnames]{xcolor}
\usepackage[T1]{fontenc}	

\usepackage{tikz}
\usetikzlibrary{arrows, positioning, calc, intersections}
\usetikzlibrary{decorations.pathreplacing, decorations.markings}

\usepackage{hyperref}
\hypersetup{
    unicode=false,          		
    pdftoolbar=true,        		
    pdfmenubar=true,        	
    pdffitwindow=false,     		
    pdfstartview={FitH},    		
    linktocpage=true,			
    pdfnewwindow=true,      	
    colorlinks=true,       			
    linkcolor=blue,          		
    citecolor=PineGreen,    	
    filecolor=magenta,      		
    urlcolor=cyan           		
}

\theoremstyle{plain}
\newtheorem{theorem}{Theorem}[section]

\newtheorem{lemma}[theorem]{Lemma}
\newtheorem{proposition}[theorem]{Proposition}

\theoremstyle{definition}

\newtheorem{problem}[theorem]{Problem}
\newtheorem{RHP}[theorem]{Riemann-Hilbert Problem}

\theoremstyle{remark}
\newtheorem{remark}[theorem]{Remark}

\numberwithin{figure}{section}
\numberwithin{equation}{section}

\allowdisplaybreaks

\usepackage{float}

\usepackage[font=small, format=plain, labelfont=bf, up, textfont=it, up]{caption}
\usepackage{subcaption}

\DeclareMathOperator{\ad}{ad}

\DeclareMathOperator{\real}{Re}
\DeclareMathOperator{\imag}{Im}

\DeclareMathOperator{\res}{res}
\DeclareMathOperator{\dist}{dist}

\newcommand{\lp}{\left(}
\newcommand{\rp}{\right)}

\DeclareFontFamily{U}{mathx}{\hyphenchar\font45}
\DeclareFontShape{U}{mathx}{m}{n}{
      <5> <6> <7> <8> <9> <10>
      <10.95> <12> <14.4> <17.28> <20.74> <24.88>
      mathx10
      }{}
\DeclareSymbolFont{mathx}{U}{mathx}{m}{n}
\DeclareFontSubstitution{U}{mathx}{m}{n}
\DeclareMathAccent{\widecheck}{0}{mathx}{"71}
\DeclareMathAccent{\wideparen}{0}{mathx}{"75}

\begin{document}


\newcommand{\purpletext}[1]{{\color{purple} #1}}
\newcommand{\bluetext}[1]{ #1 }								

\title[Derivative Nonlinear Schr\"{o}dinger Equation]{The Derivative Nonlinear Schr\"{o}dinger Equation: Global well-posedness and soliton resolution}
\author{Robert Jenkins}
\author{Jiaqi Liu}
\author{Peter Perry}
\author{Catherine Sulem}
\address[Jenkins]{Department of Mathematics, Colorado State University, Fort Collins, Colorado 80523-1874}
\address[Liu]{Department of Mathematics, University of Toronto, Toronto, Ontario M5S 2E4, Canada}
\address[Perry]{ Department of Mathematics,  University of Kentucky, Lexington, Kentucky 40506--0027}
\address[Sulem]{Department of Mathematics, University of Toronto, Toronto, Ontario M5S 2E4, Canada }
\date{\today}

\dedicatory{Dedicated to Walter Strauss, with friendship and admiration}
\begin{abstract}

We review recent results on global wellposedness and long-time behavior 
of smooth solutions to the derivative nonlinear Schr\"{o}dinger (DNLS) 
equation. Using the integrable character of DNLS, we show how the inverse 
scattering tools and the method of Zhou \cite{Zhou89-1} for treating spectral 
singularities lead to global wellposedness for general initial conditions in the weighted Sobolev space $H^{2,2}(\mathbb{R})$. For generic initial data that can support bright solitons but exclude spectral singularities, we prove the \emph{soliton resolution conjecture}: the solution is asymptotic, at large times, to a sum of localized solitons and a dispersive component, Our results also show that soliton solutions of DNLS are asymptotically stable.
\end{abstract}

\maketitle

\tableofcontents

%
%
%
%

%
%

\tikzset{->-/.style={decoration={
  markings,
  mark=at position .55 with {\arrow[scale=0.6]{triangle 45}} },postaction={decorate}}
}

%
%

%
%

%
%

\newcommand{\solfig}{
\hspace*{\stretch{.5}}
\begin{minipage}[c]{.45\textwidth}
\vspace{0pt}
\begin{tikzpicture}[scale=0.5,every node/.style={scale=0.7}]						
\coordinate (shift) at (-40:1.5);
\coordinate (x1) at (0.7,0);
\coordinate (x2) at (2.1,0);
\coordinate (topleft) at ($(x1)+(95:6)$);
\coordinate (bottomleft) at ($(x1)+(-120:6)$);
\coordinate (topright) at ($(x2)+(60:6)$);
\coordinate (bottomright) at ($(x2)+(-85:6)$);
\coordinate (C) at ($ (x1)!.5!(x2)+(77.5:.5)$);

\begin{scope}
  \clip (-4,-4) rectangle (4,4);
  \path[name path=top] (-4,4) -- (4,4);
  \path[name path=bottom] (-4,-4) -- (4,-4);
  \path[name path=right] (4,4) -- (4,-4);  	
  \path [fill=gray!10] (x1) -- (topleft) -- (topright) -- 
    (x2) -- (bottomright) -- (bottomleft) -- (x1);
  \draw[thick,name path=leftcone] (bottomleft) -- (x1) -- (topleft);
  \draw[thick,name path=rightcone] (bottomright) -- (x2) -- (topright);
    \draw [name path=axis][->] (-4,0) -- (3.3,0)
    	node[label=0:$x$] {};
    \draw [][->] (-0.5,-4) -- (-0.5,3.0)
    	node[label=180:$t$] {};
    \path [name intersections={of= leftcone and top, by=TL}];
    \path [name intersections={of= rightcone and right, by=TR}];
    \path [name intersections={of= leftcone and bottom, by=BL}];
    \path [name intersections={of= rightcone and bottom, by=BR}];
\end{scope}
    \node [fill=black, inner sep=1.5pt, circle, label=-70:$x_2$] at (x2) {};
    \node [fill=black, inner sep=1.5pt, circle, label=-177:$x_1$] at (x1) {};
	\node [label=180:${x-v_1 t = x_1}$] at ( $(TL) +(.3,-.2) $ ) {};
	\node [label=90:${x-v_2 t = x_2}$] at ( $(TR) +(-.3,-.2 ) $) {};
    \node [label=170:${x-v_2 t = x_1}$] at (BL) {};
    \node [label=10:${x-v_1 t = x_2}$] at (BR) {};
    \node at (C) {$\mathcal{S}$};
\end{tikzpicture}
\end{minipage}
\hspace*{\stretch{1}}
\begin{minipage}[c]{.45\textwidth}
\vspace{0pt}
\begin{tikzpicture}[scale=0.50, every node/.style={scale=0.7}]			
	\coordinate (v1) at (0.75,0);
	\coordinate (v2) at (-2,0);
	\path [fill=gray!10] ($(v1)+(0,-1)$) -- ($(v1)+(0,7)$) 
	--($(v2)+(0,7)$) -- ($(v2)+(0,-1)$) --  ($(v1)+(0,-1)$);
	\draw[thick] ($(v1)+(0,-.4)$) -- ($(v1)+(0,7)$);
	\draw[thick] ($(v2)+(0,-.4)$) -- ($(v2)+(0,7)$);
	\draw[][->] (-4,0) -- (4,0) 
		node[label=0:$\real \lam$] {};
	\node [label=${-v_1/4}$] at (.9,-1.25) {};
	\node [label=${-v_2/4}$] at (-2.2,-1.25) {};
\node [fill=black, inner sep = 1pt, circle, label=-90:$\lambda_1$] at (3,6) 			{};
\node [fill=black, inner sep = 1pt, circle, label=-90:$\lambda_2$] at (1.5,5) 		{};
\node [fill=black, inner sep = 1pt, circle, label=-90:$\lambda_3$] at (0,5.5) 		{};
\node [fill=black, inner sep = 1pt, circle, label=-90:$\lambda_5$] at (-3.3,4.4) 	{};
\node [fill=black, inner sep = 1pt, circle, label=-90:$\lambda_8$] at (-0.8,0.8) 	{};
\node [fill=black, inner sep = 1pt, circle, label=-90:$\lambda_4$] at (-2.5,6.5) 	{};
\node [fill=black, inner sep = 1pt, circle, label=-90:$\lambda_6$] at (-1.6,3) 		{};
\node [fill=black, inner sep = 1pt, circle, label=-90:$\lambda_9$] at (2,1) 			{};
\node [fill=black, inner sep = 1pt, circle, label=-90:$\lambda_7$] at (-3.8,1.1) 	{};
\node [fill=black, inner sep = 1pt, circle, label=-90:$\lambda_{10}$] at (3.8,2.3) 	{};
\end{tikzpicture}
\end{minipage}
\hspace*{\stretch{.5}}
}

%
%

\newcommand{\qsol}{q_{\mathrm{sol}}}
\newcommand{\usol}{u_{\mathrm{sol}}}
\newcommand{\tqsol}{\widetilde{q}_{\mathrm{sol}}}

\newcommand{\poles}{\Lambda}
\newcommand{\negint}{ I^{-}_{\xi,\sgnt}} 						
\newcommand{\sgnt}{\eta}

%
%

\newcommand{\sidenote}[1]{\marginpar{\scriptsize{\color{purple} #1}}}

\newcommand{\bff}[1]{\bm{#1}}

\newcommand{\eps}{\varepsilon}
\newcommand{\lam}{\lambda}
\newcommand{\sig}{\sigma}

\newcommand{\darr}{\downarrow}
\newcommand{\rarr}{\rightarrow}

\newcommand{\Lam}{\Lambda}

\newcommand{\dee}{\partial}
\newcommand{\dbar}{\overline{\partial}}
\newcommand{\qbar}{\overline{q}}

\newcommand{\dotarg}{\, \cdot \, }

\newcommand{\ba}{\breve{a}}
\newcommand{\balpha}{\breve{\alpha}}
\newcommand{\bb}{\breve{b}}
\newcommand{\bbeta}{\breve{\beta}}
\newcommand{\br}{\breve{r}}
\newcommand{\bbC}{\mathbb{C}}

\newcommand{\zetabar}{\overline{\zeta}}

\newcommand{\scrB}{\mathscr{B}}

\newcommand{\C}{\mathbb{C}}
\newcommand{\I}{\mathbb{I}}
\newcommand{\N}{\mathbb{N}}
\newcommand{\R}{\mathbb{R}}
\newcommand{\Z}{\mathbb{Z}}

\newcommand{\calA}{\mathcal{A}}
\newcommand{\calB}{\mathcal{B}}
\newcommand{\calD}{\mathcal{D}}
\newcommand{\calC}{\mathcal{C}}
\newcommand{\calG}{\mathcal{G}}
\newcommand{\calI}{\mathcal{I}}
\newcommand{\calL}{\mathcal{L}}
\newcommand{\calJ}{\mathcal{J}}
\newcommand{\calM}{\mathcal{M}}
\newcommand{\calQ}{\mathcal{Q}}
\newcommand{\calR}{\mathcal{R}}
\newcommand{\calS}{\mathcal{S}}
\newcommand{\calT}{\mathcal{T}}
\newcommand{\calU}{\mathcal{U}}

\newcommand{\One}{\mathbf{1}}

\newcommand{\kbar}{\overline{k}}
\newcommand{\ubar}{\overline{u}}
\newcommand{\zbar}{\overline{z}}

\newcommand{\bfA}{\mathbf{A}}
\newcommand{\bbR}{\mathbb{R}}
\newcommand{\bfr}{\mathbf{r}}
\newcommand{\bfM}{\mathbf{M}}
\newcommand{\bfN}{\mathbf{N}}
\newcommand{\bft}{\mathbf{t}}

\newcommand{\tildJ}{\widetilde{J}}
\newcommand{\brevcalJ}{\breve{\calJ}}

\newcommand{\dint}{\displaystyle{\int}}

\newcommand{\norm}[2][ ]{\left\Vert #2 \right\Vert_{#1}}
\newcommand{\bigO}[2][ ]{\mathcal{O}_{#1}\left( #2 \right)}

\newcommand{\upmat}[1]{\ttwomat{0}{#1}{0}{0}}
\newcommand{\lowmat}[1]{\ttwomat{0}{0}{#1}{0}}

\newcommand{\twomat}[4]
{
\begin{pmatrix}
#1	&	#2	\\
#3	&	#4
\end{pmatrix}
}

\newcommand{\Twomat}[4]
{
\begin{pmatrix}
#1	&	#2	\\[5pt]
#3	&	#4
\end{pmatrix}
}


\newcommand{\ttwomat}[4]
{
\begin{pmatrix}
#1	&	#2	\\[3pt]
#3	&	#4
\end{pmatrix}
}

\newcommand{\diagmat}[2]
{
\begin{pmatrix}
#1		&	0	\\
0		&	#2
\end{pmatrix}
}

\newcommand{\offmat}[2]
{
\begin{pmatrix}
0	&	#1		\\
#2	&	0
\end{pmatrix}
}

\newcommand{\twovec}[2]
{
	\left(
		\begin{array}{c}
			#1		\\
			#2
		\end{array}
	\right)
}

\newcommand{\Twovec}[2]
{
	\left(
		\begin{array}{c}
			#1		\\
			\\
			#2
		\end{array}
	\right)
}

%
%

\newcommand{\sol}{{\mathrm{sol}}}

%
%


\section{Introduction}
This paper is devoted to a review of   recent results on global wellposedness and large-time asymptotics of the 
the derivative nonlinear Schr\"odinger   (DNLS)  equation
\begin{equation}
\label{DNLS}
iu_t + u_{xx} - i \eps \left( |u|^2 u \right)_x = 0, \quad x\in \bbR
\end{equation}
with initial condition
\begin{equation} \label{IC}
u(x,0)=u_0(x)
\end{equation}
and  $\varepsilon = \pm 1$. The transformation $u(x,t) \mapsto u(-x,t)$ maps solutions of \eqref{DNLS} with $\eps=-1$ to solutions of \eqref{DNLS} with $\eps=1$.


The DNLS equation is a canonical  dispersive equation that can be obtained  in a  long-wave, weakly nonlinear scaling regime
from the  one-dimensional  compressible magneto-hydrodynamic (MHD) equations  in the presence of the Hall effect  (Mj\o lhus \cite{M76}, 
see also Champeaux et.\  al.\  \cite{CLPS99}). This derivation is  somehow  similar to that of the Korteweg de Vries equation from 
the water wave problem and we will discuss in more details  later.

Under gauge transformations,  the DNLS equation takes  equivalent forms that can be useful in various contexts and purposes (see, for example, Wadati-Sogo \cite{WS83} and references therein).
The form that we will adopt here is 
\begin{equation}
\label{DNLS2}
iq_t + q_{xx} + i\eps q^2 \qbar_x + \frac{1}{2}|q|^4 q = 0
\end{equation}
obtained under the invertible gauge transformation
\begin{equation}
\label{DNLS:q.gauge}
q(x,t) = u(x,t) \exp\left(-i\eps \int_x^{\infty} |u(y,t|^2 \, dy \right).
\end{equation}
In the context of inverse scattering, the latter equation has a Lax representation whose direct scattering problem is more easily normalized at spatial infinity.

DNLS is invariant under the scaling transformation
\begin{equation}\label{scaling}
u \rightarrow u_\lambda = \lambda^{-1/2} u({x/\lambda}, {t}/{\lambda^2}).
\end{equation}
In particular, it is  $L^2$-critical in the sense that $\| u_\lambda \|_{L^2} = \| u \|_{L^2}$.

Local wellposedness of smooth solutions (in the Sobolev spaces $H^s$, $s>3/2$)
was established by Tsutsumi and Fukuda \cite{TF80} and later extended to solutions with low regularity  (in the Sobolev space  $H^{1/2})$
by Takaoka \cite{T99}. 
For initial conditions in the energy space $H^1$, Hayashi-Ozawa \cite{HO92} proved that solutions exist in $H^1$ for all time if $\|u_0\|_{L^2} < \sqrt{2\pi}$.
The proof uses conservation laws and the optimal constant in the Gagliardo-Nirenberg inequality in a setting similar to that of  Weinstein for 
$L^2$-critical NLS equations. Colliander-Keel-Staffilani-Takaoka-Tao \cite{CKSTT02}  extended this result to 
$u_0\in H^{1/2+\epsilon}$.  
More recently,  Wu \cite{W15}  and Guo-Wu \cite{GW17} increased  the upper bound $2\pi$ to $ 4\pi$, with  initial conditions  respectively in $u_0\in H^{1}$
and  $H^{1/2}$. Our results give global existence in $H^{2,2}(\R)$ with no restriction on the $H^{2,2}(\R)$ norm. obtaining a ``global'' result at the cost of assuming greater regularity and decay for the initial data.

In term of minimal regularity, this result is  optimal. Indeed, for $u_0 \in H^s$, $s < 1/2$,  the Cauchy problem is ill-posed, in the sense 
that uniform continuity with respect to the initial conditions fails \cite{T99}.

It is well-known (Weinstein \cite{We83}) that for focusing $L^2$-critical NLS equations, the optimal  (strict) upper bound on initial $L^2$-norm to ensure global
well-posedness is given by the $L^2$-norm of the ground state (i.e. the unique positive solution of $\Delta R-R+R^{4/n}=0$, $n$ being the space dimension). Here,  It is also of interest to relate the above constants $\sqrt{2\pi}$ and $\sqrt{4\pi}$
to  the $L^2$-norm of solitary waves. 

Fixing $\varepsilon = - 1$  in eq. \eqref{DNLS}, the DNLS equation has a 
two-parameter family of solitary waves, in the form  ($\omega > c^2/4$) 
\begin{equation}\label{u 1sol}
  u_{\omega, c}(x,t) = \varphi_{\omega, c}(x-ct)\exp i \left\{ \omega t
    +\frac{c}{2}(x-ct)- \frac{3\eps}{4}\int^{x-ct}_{-\infty}\varphi^{2}(\eta)d\eta\right\}
\end{equation}
where 
$$ \varphi_{\omega, c}(y) = \sqrt{ \frac{(4\omega-c^2)}{\omega^{1/2} (\cosh(\sqrt{4\omega-c^2}y)-\frac{c}{2\sqrt{\omega}})}}$$
is the  unique  solution of 
$$  - \partial_{y}^2\varphi + \left(\omega -
  \frac{c^2}{4}\right)\varphi +\frac{c}{2}|\varphi |^{2}\varphi - \frac{3}{16}|\varphi|^{4}\varphi= 0.$$
  These solutions decay exponentially fast at infinity and are called bright solitons.
In the limiting case  $ \sqrt{\omega} = c/2$,   the profile $\varphi_{\omega, c}$ reduces to 
        $$ \varphi_{\omega, 2 \sqrt{\omega}}(x) = \frac{2\omega\sqrt{2}}{\sqrt{1+4\omega x^2},}$$
        usually called  the  lump or  algebraic soliton.
          The $L^2$-norm of $u_{\omega,c}$ can be calculated explicitely:
 $$ \|u_{\omega,c}\|_{L^2}^2 =   8 \tan^{-1} { \sqrt{\frac{2 \sqrt{\omega}+ c}{ 2 \sqrt{\omega} -c }}}     \le 4\pi $$
 where the limiting value $4\pi$ corresponds to the square of the 
   $L^2$-norm of algebraic soliton.

The  orbital stability of   bright 1-soliton  solution, that is the stability in $H^1$ up to the two transformations, translation and multiplication by a constant phase that leave the equation invariant, was established  by Colin-Ohta \cite{CO06} and extended to $N$-soliton solutions by  Le Coz-Wu \cite{LW18}. 
In the case of the 
algebraic soliton, Kwong-Wu  \cite{KW18} proved that  
orbital stability  holds in the following sense: 
If the initial condition $u_0$ is close in $H^1$-norm
to the algebraic soliton $R= u_{\omega, 2 \sqrt{\omega}}$, then, for $t\in [0,t_*)$, where $t_*$ is the maximal time of existence,
there exist $\theta(t)$, $y(t)$ and $\lambda(t)$ such that the solution $u(\cdot,t)$ remains close to $e^{i\theta(t)} R_{\lambda(t)} (\cdot-y(t))$ in $H^1$ norm,
where $R_\lambda$ is obtained from $R$ by the scaling transformation \eqref{scaling}.

We now recall known results concerning the long-time behavior of DNLS solutions. The first results go back to 
Hayashi, Naumkin and Uchida (1999) where the authors consider a class of  one-dimensional 
nonlinear Schr\"odinger  equations with 
general nonlinearities containing first-order derivatives.  They prove a global existence result for  smooth initial conditions that are small in some weighted  Sobolev spaces, as well as  a  time-decay rate. Their analysis provides the existence of asymptotic states $u^\pm \in L^2 \cap L^\infty$ and 
real valued functions $ {g}^\pm \in L^\infty$ such that

 
  \ 
\begin{equation*}
u (x,t)   \sim \frac{1}{\sqrt{t}} u^\pm(x/(2t)) ~\exp \left( \frac{ix^2}{4t}  \pm i {g}^\pm(z/(2t)) \log |t| \right)  +O( |t|^{-1/2-\alpha})
\end{equation*}
uniformly in $x \in \bbR$, with  $0 < \alpha < 1/4$.
 Kitaev and Vartanian \cite{KV97,KV99} obtained large-time asymptotic expansions using inverse scattering methods, but imposed a small-norm assumption on the scattering data which we were able to avoid. Their work provided a number of valuable hints for ours.

The goal of this review article is to present new results obtained during the last five years in a series of papers \cite{JLPS19,JLPS18a,JLPS18b,LPS16}
and the thesis \cite{Liu17} of the second author that  have provided
answers to fundamental questions of  global existence of solutions for large $L^2$ data,  long-time behavior of solutions,
asymptotic stability of solitons and more generally the soliton resolution conjecture. 
Our approach is based on a central   structural property of DNLS, discovered by  Kaup and Newell \cite{KN78}  that it is 
solvable through the inverse scattering method. In this pioneering work, the authors 
establish the main elements of the inverse scattering analysis.  In particular, they find the Lax pair, analyze the linear spectral flow and derive the soliton solutions.

Before setting the tools of inverse scattering and stating our results, we briefly present the formal derivation of  DNLS from the original
Hall-MHD equations  (Mjolhus \cite{M76}).

\subsection{DNLS as a long-wave, small-amplitude model}
The starting point is the compressible MHD equations in the presence of the Hall term in the Ohm law (responsible for the dispersive character of the equations),  assuming  a uniform background magnetic field along
the $x$-axis and variations of the different fields in the $x$ variable only. 

Differently from  the long-wave regime  of an ideal potential flow, the  long-wave asymptotics does not  lead to the KdV equation
but to two coupled equations for the two components of the magnetic field transverse to the propagation conveniently written, using
complex notation in one single equation, the so-called derivative NLS equation. In this section, we summarize the main steps of the derivation
following Champeaux et al \cite{CLPS99}.

Denoting $\rho$ the density of the fluid, $u= (u_1, u_2, u_3)$ its velocity and $B = (b_1, b_2, b_3)$ the magnetic field,  we introduce the complex transverse  variables $v= u_2- \sigma u_3$,  $b= b_2 -i\sigma b_3$ where $\sigma= \pm 1$ corresponds to the right-hand or left-hand 
circularly polarized waves.  The longitudinal component $b_1$ can be absorbed into the background magnetic field  and is taken to be equal  to 1.  
The system takes the form
\begin{align*}
\partial_t \rho + \partial_x(\rho u_1) &=0 \\
 \rho(\partial_t u_1 + u_1\partial_x u_1 ) &=  - \left( \frac{\beta}{\gamma} \rho^\gamma +\frac{1}{2} |b|^2 \right)\\
 \rho(\partial_t v + u_1\partial_x v ) &=\partial_x b \\
\partial_t b + u_1\partial_x  ( u_1 b-v )  &=i  \frac{\sigma}{R} \partial_x \left( \frac{1}{\rho} \partial_x b\right)
\end{align*}
where  $R, \beta, \gamma$ are constants. $R$ is the ion-gyrofrequency,   $\beta,  \ne 1$ is   the square of the ratio  of the sonic and the Alfv\'en speeds  and  $\gamma$ the   polytropic gas constant.
This is a dispersive system and the dispersion relation is (assuming $k>0$)
$$ \omega(k) = \frac{\sigma}{2R} k^2 + k\sqrt{ 1 + \frac{k^2}{4R^2}}.$$

To capture the long-wave, small amplitude regime, we introduce the rescaled independent variables 
$$ \xi =\varepsilon( x-t), \quad \tau = \varepsilon t$$
and expand the dependent variables in the form
\begin{align*} 
& \rho= 1+\varepsilon \rho_1 +\cdots,  &u_1 = \varepsilon (u_1^{(1)} +\varepsilon u_1^{(2)} +\cdots ), \cr
& v= \varepsilon^{1/2} (v_1 + \varepsilon v_2 +\cdots),   &b= \varepsilon^{1/2} (b_1 + \varepsilon b_2 +\cdots).
\end{align*}
After substitution in the MHD equations and identification, one gets at order $O(\varepsilon^{3/2})$
$$\partial_\xi(b_1+v_1) =0,$$
thus $ b_1 =-v_1$. At order $O(\varepsilon^{3/2})$, 
$$u_1^{(1)} = \rho_1 = \frac{|b_1|^2}{2(1-\beta)}.$$
At order $O(\varepsilon^{5/2})$, the resulting equation for $b_1$ is (dropping the index 1)
$$\partial_\tau b +\frac{i}{2R} \partial_{\xi\xi} b + \frac{1}{4 ( 1-\beta) } \partial_\xi (|b|^2b)=0,$$
whixh is the DNLS equation.

A central problem in nonlinear PDEs 
has been to describe in which sense certain canonical asymptotic equations provide approximations
to physical problems, and to study rigorously the behavior of solutions in these asymptotic
limits. Justification of long-wave models such as Boussinesq and KdV equations derived from the
original Euler equations for water waves has been the object of intense studies in the last 40 years  \cite{BLS08,C85,L13,SW00} and
are well-understood.
The validity of DNLS has been tested against direct numerical simulations of the Hall-MHD equations  in \cite{CLPS99} but a rigorous analysis of the asymptotics has not been done  so far.

\subsection{DNLS as an integrable system}

Kaup and Newell \cite{KN78} showed that equation \eqref{DNLS} is the consistency condition for the overdetermined system
\begin{align}
\label{DNLS:x}
\bff{\psi}_x	&=	\left(-i \zeta^2 \sigma_3 +\zeta \bff{U} \right) \bff{\psi} \\[5pt]
\label{DNLS:t}
\bff{\psi}_t		&=	\left( \zeta^4 \sigma_3 +2i\zeta^3 \bff{U} + \zeta^2 \eps |u|^2 \sigma_3 +   i\zeta \eps |u|^2 \bff{U} - i \zeta \bff{U}_x \right) \bff{\psi}.
\end{align}
Here $\bff{\psi}(x,t)$ is an unknown $2 \times 2$ matrix-valued function of $(x,t)$ and $\sigma_3$ is the Pauli matrix
$$ 
\sigma_3 = 
\begin{pmatrix}
1 	&	0	\\ 0	& -1 
\end{pmatrix}.
$$ 
The matrix $\bff{U}=\bff{U}(x,t)$ depends on the putative solution $u(x,t)$ via the formula
$$ \bff{U}(x,t) = 
		\begin{pmatrix}
			0	&	u(x,t)	\\
			\eps \overline{u(x,t)}	&	0
		\end{pmatrix}.
$$
The system \eqref{DNLS:x}--\eqref{DNLS:t} is called a \emph{zero-curvature representation} for the DNLS equation; one obtains \eqref{DNLS} by cross-differentiating in $x$ and $t$ and appealing to Clairaut's Theorem. 
This special symmetry means that \eqref{DNLS} is an isospectral\footnote{To be precise, the problem \eqref{DNLS:x} is not a spectral problem but rather an \emph{operator pencil} since the right-hand side depends polynomially on the spectral parameter. This turns out to cause some technical complications but is otherwise an inessential matter.} flow for the problem \eqref{DNLS:x} whose spectral and scattering evolve linearly with time in a manner dictated by \eqref{DNLS:t}. 

In his 1983 thesis \cite{Lee83} and subsequent paper \cite{Lee89}, J.-H.\ Lee applied the Beals-Coifman \cite{BC84} approach to inverse scattering to study the initial value problem for \eqref{DNLS} for generic initial data in the Schwartz class of smooth functions of rapid decrease. Following the lead of Kaup and Newell \cite{KN78}, Lee actually considered the  gauge-equivalent evolution equation\footnote{For a discussion of the gauge equivalence see, for example, \cite[Appendix A]{LPS16}.} \eqref{DNLS2}
and its zero curvature representation
\begin{align}
\label{DNLS2:x}
\bff{\psi}_x	&=	\left(-i\zeta^2 \sigma_3 + \zeta \bff{Q} + \bff{P} \right) \bff{\psi} \\
\label{DNLS2:t}
\bff{\psi}_t		&=	\Bigl( 2\zeta^4 \sigma_3 + \zeta^2 \eps|q|^2 \sigma_3  + 2i\zeta^2 \bff{Q} - \zeta \sigma_3 \bff{Q}_x   \Bigr. \\
\nonumber
					&	\quad \qquad \Bigl. + \frac{i}{2} \sigma_3( \eps q_x \qbar - q \qbar_x) - \frac{1}{4} \sigma_3 |q|^4 \Bigr) \bff{\psi}
\end{align}
where
$$
\bff{Q}(x,t) = 
\begin{pmatrix}
0								&	q(x,t) 	\\
\eps \overline{q(x,t)}		&	0
\end{pmatrix},
\quad
\bff{P}(x,t)	=	-\eps\frac{i}{2} |q(x,t)|^2 \sigma_3. 
$$
Equation \eqref{DNLS2:x} defines a spectral problem as follows; for this discussion, we will write $q(x)$ for $q(x,t)$, ignoring the dependence on the $t$ variable.  The reader familiar with the ZS-AKNS system associated to the cubic nonlinear Schr\"{o}dinger equation (see \cite{AKNS74,ZS72}) will note that the scattering and inverse scattering theory for \eqref{DNLS2:x} are very similar, although the quadratic dependence on the spectral parameter in \eqref{DNLS2:x} introduces some complications.

First, note that if
 $\bff{P}=\bff{Q}=\bff{0}$,
\eqref{DNLS2:x} admits bounded solutions of the form $\exp(-ix\zeta^2 \sigma_3)$ so long as 
$\imag \zeta^2 = 0$, that is, for $\zeta$ in the cross-shaped region shown in Figure \ref{fig:Sigma}. If $q \in L^1(\R) \cap L^2(\R)$, there exist unique solutions $\bff{\psi}^\pm(x,\zeta)$ 
of \eqref{DNLS2:x} with $$\lim_{x \to \pm \infty}  \bff{\psi}(x,\zeta) e^{i\zeta^2 x \sigma_3} = \I, \quad \I = 
\begin{pmatrix}
1	&	0	\\	
0	&	1
\end{pmatrix}.
$$

\begin{center}
\begin{figure}[H]
\caption{The Contour $\Sigma$ and the regions $\Omega^\pm$}

\medskip 

\begin{tikzpicture}
\draw[white,fill=gray!50]	(0,0)	rectangle	(-3,3);
\draw[white,fill=gray!50]	(0,0)	rectangle	(3,-3);
\node at		(-1.5,1.5)		{$\Omega^-$};
\node at		(1.5,-1.5)		{$\Omega^-$};
\node at		(1.5,1.5)		{$\Omega^+$};
\node at		(-1.5,-1.5)		{$\Omega^-$};
\draw[thick,->-]	(0,0)		--	(-3,0);
\draw[thick,->-]	(0,0)		--	(3,0);
\draw[thick,->-]	(0,3)		--	(0,0);
\draw[thick,->-]	(0,-3)		--	(0,0);
\node[right] 	at 	(3,0)	{$\real \zeta$};
\node[above]	at	(0,3)	{$\imag \zeta$};
\end{tikzpicture}
\label{fig:Sigma}
\end{figure}
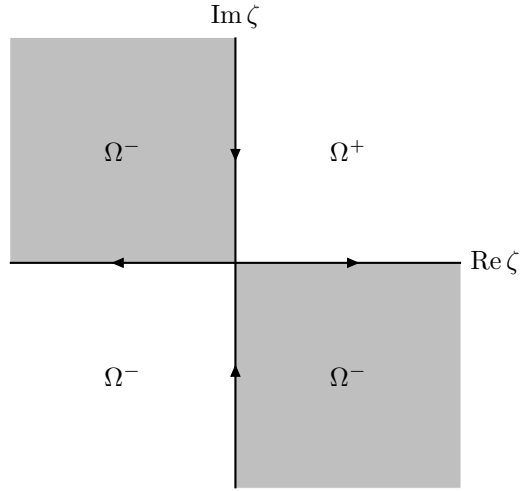
\end{center}
The solutions $\bff{\psi}^\pm(x,\zeta)$ are called \emph{Jost solutions}. Since \eqref{DNLS2:x} takes the form $\bff{\psi}_x = \bff{B}(x)\psi$ where $\bff{B}$ is a traceless matrix, it follows that $\det \bff{\psi}^\pm(x,\zeta)$ is constant, hence $\det \bff{\psi}^\pm(x,\zeta) = 1$. It also follows from the form of the equation that any two nonsingular solutions $\bff{\psi}_1$ and $\bff{\psi}_2$ of \eqref{DNLS2:x} satisfy $\bff{\psi}_1 = \bff{\psi}_2 \bff{A}$ for a nonsingular matrix $\bff{A}$. Thus, there is a matrix $\bff{T}(\zeta)$ of determinant one so that
\begin{equation}
\label{DNLS2:T}
\bff{\psi}^+(x,\zeta) = \bff{\psi}^-(x,\zeta) \bff{T}(\zeta), \quad 
\bff{T}(\zeta) = \begin{pmatrix}
a(\zeta)	&	\bb(\zeta)	\\
b(\zeta)	&	\ba(\zeta)
\end{pmatrix}.
\end{equation}
The entries of $\bff{T}$ obey the symmetries
\begin{equation}
\label{symmetry}
 a(\zeta) = \overline{\ba(\zetabar)}, \quad b(\zeta) = \eps \overline{\bb(\zetabar)}, \quad
a(-\zeta) = a(\zeta), \quad b(\zeta) = -b(-\zeta)
\end{equation}
together the determinant relation $a(\zeta) \ba(\zeta) - b(\zeta) \ba (\zeta) = 1$.  
 Assuming that $a$ and $\ba$ are zero-free in their respective regions of definition, one can define the reflection coefficients
\begin{equation}
\label{rba}
r(\zeta) = \bb(\zeta)/a(\zeta), \quad \br(\zeta) = b(\zeta)/\ba(\zeta). 
\end{equation}
It can be shown that $a(\zeta)$ extends to an analytic function on the region $\Omega^-$ (see Figure \ref{fig:Sigma}) while $\ba(\zeta)$ extends to an analytic function on $\Omega^+$. Zeros of $a$ and $\ba$ are associated to soliton solutions of \eqref{DNLS2}.
Roughly speaking, the \emph{direct scattering map} is the mapping from given $q$ to the functions $a$, $\ba$, $b$, $\ba$. We will give a more precise formulation later.

To define the inverse map, we need to introduce the notion of Beals-Coifman solutions, following \cite{BC84,Lee83}. If we factor a solution $\psi$ of \eqref{DNLS2:x} as
$$\bff{\psi}(x,\zeta) = \bff{M}(x,\zeta) e^{-i\zeta^2 x \sigma_3},$$ then $M$ solves the problem
\begin{equation}
\label{DNLS2:M}
\frac{d}{dx} \bff{M}(x,\zeta) = - i\zeta^2 \ad \sigma_3 (\bff{M}) + \zeta \bff{Q}(x) \bff{M} + \bff{P}(x) \bff{M}.
\end{equation}
A solution of \eqref{DNLS2:M} for $\zeta$ with $\imag \zeta^2 \neq 0$ is called a left (resp.\ right) Beals-Coifman solution if $\lim_{x \to -\infty} \bff{M}(x,\zeta) = \I$ (resp.\ $\lim_{x \to \infty} \bff{M}(x,\zeta) = \I$) and $\bff{M}(x,\zeta)$ is bounded as $x \to \infty$ (resp.\ $\bff{M}(x,\zeta)$ is bounded as $x \to -\infty$). It turns out that the Beals-Coifman solutions exist and are unique so long as $\zeta$ is not a zero of $a(\zeta)$ or $\ba(\zeta)$. 

If we now re-instate the dependence of $q(x,t)$ (and therefore the scattering data $a(\zeta,t)$, $b(\zeta,t)$ etc.) on $t$, we can use the second equation \eqref{DNLS2:t} of the zero curvature representation to compute the time-dependence of the scattering data $a(\zeta,t)$ and $b(\zeta,t)$ if $q(x,t)$ is a solution. One finds that
\begin{equation}
\label{DNLS2:ab-ev}
\dot{a}(\zeta,t) = 0, \quad \dot{b}(\zeta,t)= -4i\zeta^4 b(\zeta,t). 
\end{equation}
Thus, the scattering data obey a \emph{linear} evolution. Given a way of recovering $q(x,t)$ from $a(\zeta,t)$ and $b(\zeta,t)$, we can solve \eqref{DNLS2}. The recovery is implemented by the \emph{inverse scattering map} determined by a Riemann-Hilbert problem.

To describe the inverse scattering map, we once again suppress dependence of all quantities on $t$ and, in effect, describe the map acting on scattering data at a fixed time. We also make a simplifying assumption which we will carry through for our discussion of soliton resolution but drop for our discussion of global well-posedness. 
A zero of $a$ or $\ba$ is called a \emph{spectral singularity} if it lies on the real or imaginary axis, and a \emph{resonance} otherwise.  
We assume that $q(x)$ is so chosen that (i)  $a(\zeta)$ has at most finitely many zeros, 
and (ii) none of these zeros lie on the contour $\Sigma$.  
In \cite{JLPS18b}  it is shown that conditions (i) and (ii) hold for generic functions $q$.\footnote{More precisely, we show that the set of $q$ satisfying these conditions is an open and dense subset of the Sobolev space 
$$H^{2,2}(\R) = \{ q \in L^2(\R): x^2q, \,\, q'' \in L^2(\R) \}. $$} We denote by $Z$ the finite set of $\zeta \in \C \setminus \Sigma$ for which $a(\zeta)=0$ or $\ba(\zeta)=0$; it follows from the symmetries of $a$ and $\ba$ that $Z$ is a union of `quartets' of the form $(\zeta_j, -\zeta_j, \zetabar_j, -\zetabar_j)$ with $\real \zeta_j, \imag \zeta_j > 0$, as shown in Figure \ref{fig:quartet}.

\begin{center}
\begin{figure}[H]
\caption{Zeros of $a$ and $\ba$}

\medskip

\begin{tikzpicture}
\draw[very thick,->-]	(0,0)		--	(-3,0);
\draw[very thick,->-]	(0,0)		--	(3,0);
\draw[very thick,->-]	(0,3)		--	(0,0);
\draw[very thick,->-]	(0,-3)		--	(0,0);
\node[right] 	at 	(3,0)	{$\real \zeta$};
\node[above]	at	(0,3)	{$\imag \zeta$};
\draw[black,fill=black]	(1.5,1.5)	circle(0.05cm)	node[anchor=west] {$\zeta_j$};
\draw[black,fill=black]	(1.5,-1.5)	circle(0.05cm)	node[anchor=west]	{$\zetabar_j$};
\draw[black,fill=black]	(-1.5,1.5)	circle(0.05cm)	node[anchor=east]	{$-\zetabar_j$};
\draw[black,fill=black]	(-1.5,-1.5)	circle(0.05cm)	node[anchor=east]	{$-\zeta_j$};
\end{tikzpicture}
\label{fig:quartet}
\end{figure}
\end{center}

The left Beals-Coifman solutions satisfy the following Riemann-Hilbert problem (RHP). To emphaize the role of $x$ (and later $x,t$) as a parameter, we write $\bff{M}(\zeta;x)$ for $\bff{M}(x,\zeta)$.

\begin{RHP}
\label{DNLS2:RHP.left}
Given $x \in \R$, find a function $\bff{M}(\cdot;x)$ analytic on $\C \setminus (\Sigma \cup Z)$, with the following properties:
\begin{enumerate}[(i)]
\item	$\bff{M}(\zeta;x)$ has continuous boundary values $\bff{M}^\pm$ on $\Sigma$ as $\pm \imag \zeta^2 > 0$ and
$$
 \bff{M}_+(\zeta;x)  =\bff{M}_-(\zeta;x) e^{-i\zeta^2 x \ad \sigma_3} \bff{J}(\zeta), 
$$
where
$$
\bff{J}(\zeta) = 
\begin{pmatrix}
	1-b(\zeta)\bb(\zeta)/a(\zeta) \ba(\zeta)	&	\bb(\zeta)/a(\zeta)	\\[5pt]
	-b(\zeta)/\ba(\zeta)	&	1
\end{pmatrix}
$$
\item	The residue conditions
		$$ \res_{\zeta= \zeta_*} \bff{M}(\zeta;x) = \lim_{\zeta \to \zeta_*} \bff{M}(\zeta;x) \bff{J}(\zeta_*) $$
		hold, where 
		$$ 
		\bff{J}(\zeta_*) = \begin{pmatrix}
			0	&	0\\
			c_j & 0
		\end{pmatrix}, \quad \imag \zeta_*^2 >0 , \qquad
		 \bff{J}(\zeta_*) = \begin{pmatrix}
			0	&	\overline{c_j}	\\
			0	&	0
		\end{pmatrix},
		\quad \imag \zeta_*^2 < 0. 
		$$
\item	$\lim_{z \to \infty} \bff{M}(x,z) = \I$ uniformly in sectors properly contained in $\C \setminus \Sigma$. 
\end{enumerate}
\end{RHP}

 The right Beals-Coifman solutions satisfy a similar Riemann-Hilbert problem but with different jump matrices.

 To re-introduce time dependence, we replace the jump  matrix $e^{-i\zeta^2 x \ad \sigma_3} \bff{J}(\zeta)$ above by
 $e^{-i(\zeta^2 x +2 \zeta^4 t)  \ad \sigma_3} \bff{J}(\zeta)$
due to the simple time evolution \eqref{DNLS2:ab-ev} of the coefficients $a$ and $b$.
Under suitable decay and regularity assumptions for the coefficients $a$, $\ba$, $b$, $\bb$, it can be shown that the Riemann-Hilbert problem has a unique solution admitting a 
large-$\zeta$ asymptotic expansion 
$$ \bff{M}(\zeta;x,t) \sim \I + \frac{\bff{M}_1(x,t)}{\zeta} + o\left(\frac{1}{\zeta}\right). $$
By substituting this large-$\zeta$ expansion into the differential equation \eqref{DNLS2:M} satisfied by the Beals-Coifman solutions, we can read off the reconstruction formula 
\begin{equation}
\label{DNLS2:recon}
q(x,t) = \lim_{\zeta \to \infty} 2i\zeta M_{12}(\zeta;x,t). 
\end{equation}
RHP \ref{DNLS2:RHP.left} and the reconstruction formula \eqref{DNLS2:recon} define the inverse scattering map. 

We will use RHP \ref{DNLS2:RHP.left} and the Deift-Zhou steepest descent method to obtain soliton resolution fin case $a$ and $\ba$ have only finitely many zeros in $\C \setminus \Sigma$ and no zeros on $\Sigma$. 

As we will discuss, we can use Zhou's method \cite{Zhou89-2} to recast the Riemann-Hilbert problem in a form that eliminates the need for the genericity hypothesis, at the cost of obtaining rather poor estimates of large-time behavior (and in particular no results on soliton resolutions).  We can nonetheless use this formulation to prove global well-posedness with no spectral assumptions.

\subsection{An important change of variables}
To analyze the direct map (from the given potential $q_0$ to the scattering data) and the inverse map (from the scattering data  to the recovered potential) it is helpful to exploit the symmetry reduction of the spectral problem \eqref{DNLS2:x} to the spectral variable $\lam = \zeta^2$. Under the map $\zeta \mapsto \zeta^2$, the contour $\R \cup i\R$ maps to $\R$ with its usual orientation and $\Omega^\pm$ map to $\C^\pm$, and RHP
\ref{DNLS2:RHP.left} reduces to an RHP with contour $\bbR$.
We set 
$$\alpha(\lambda) = a(\zeta), \quad \beta(\lambda) = \zeta^{-1} \bb(\zeta) $$
which satisfy the relation
$  |\alpha(\lambda)|^2 +\lam |\beta(\lambda)|^2 = 1 $
as well as 
\begin{align*}
\rho(\lambda)&= \zeta^{-1}r(\zeta), \quad C_j =2c_j
\end{align*}
It can be shown that the diagonal of $\bff{M}(\zeta;x)$ is even under the reflection $\zeta \mapsto -\zeta$, while the off-diagonal is odd. Hence, 
for 
$$ 
\bff{M}^\sharp(\zeta^2;x) = \Twomat{M_{11}(\zeta;x)}
											{\zeta^{-1} M_{12}(\zeta;x)}
											{\zeta M_{21}(\zeta;x)}
											{M_{22}(\zeta;x)}
$$ 
condition (i) in  RHP
\ref{DNLS2:RHP.left}  becomes
\begin{align} \label{jump}
\bff{M}^\sharp_+(\lam;x)	&=	\bff{M}^\sharp_-(\lam;x) e^{-i\lambda x \ad(\sigma_3)} \bff{J}(\lambda)\cr
\bff{J}(\lam)				&=	\Twomat{1 - \eps \lambda |\rho(\lambda)|^2}
											{\rho(\lambda)}
											{- \eps \lambda \overline{\rho(\lambda)}}
											{1}.
\end{align}
However, this RHP is not properly normalized. A careful computation shows that 
$$ \bff{M}^\sharp(\zeta^2;x) \rightarrow \Twomat{1}{0}{(-i/2) \overline{q} }{1}$$ 
as $\zeta^2 \rightarrow \infty$ and it is more effective to consider the row-wise RHP for 
$$ \textbf{N}(z;x) = \left( N_{11}(z;x), N_{12}(z;x) \right). $$

We arrive at the following new Riemann-Hilbert problem.

\begin{RHP}
\label{RHP2.row}
Given $x,t \in \R$, a function $\rho(\lambda)$ for $\lambda \in \R$, and $\{ (\lam_j, C_j)\}_{j=1}^N$ in $(\C^{+} \times \C^\times)^N$, find a row vector-valued function 
$$\bfN(z; x,t): \C \setminus \left(\R \cup \{{\lambda_j}\}_{j=1}^N \cup \{{\overline{\lambda}_j}\}_{j=1}^N   \right)  \to \C^2$$ with the following properties:
\begin{itemize}
\item[(i)]		$\textbf{N}(z;x,t) = \begin{pmatrix} 1 & 0 \end{pmatrix} + \bigO{\dfrac{1}{\lam}}$ as $|z| \to \infty$,
\item[(ii)]  $\textbf{N}$ has continuous boundary values $\textbf{N}_\pm$ for $\lam \in \R$ and
				$$ \textbf{N}_+(\lam;x,t) = \textbf{N}_-(\lam;x,t) e^{-it\theta(\lam,x/t) \ad (\sigma_3)} \bff{J}(\lam),
				\quad 
				\bff{J}(\lam) = \begin{pmatrix}
								1-\eps\lam |\rho(\lam)|^2 	&	 \rho(\lam)	\\
								-\eps \lam \overline{\rho(\lam)}	&	1
							\end{pmatrix},
				$$
				where
				\[
				\theta =\theta(\lambda, \xi)  = 2 \lambda^2 + \lam \xi
				\] 
\item[(iii)]	For each $\lam$ equal to $\lam_j$ or $\overline{\lambda}_j$, 
				$$ \res_{z = \lam} \textbf{N}(\lam;x,t) = \lim_{z \to \lam} \textbf{N}(z;x,t) e^{-it\theta(z,x/t) \ad (\sigma_3)} \bff{J}(\lam) $$
				where for each $\lambda= \lam_j$
				$$ \bff{J}(\lam_j)  = \lowmat{\lam_j C_j}, \quad \bff{J}( \overline{\lam_j}) = \upmat{\eps \overline{C_j}}. $$
\end{itemize} 
\end{RHP}
 A consequence of the change of variables which is central to the long-time behavior analysis  is that   $\theta$ has only one stationary point.  Given the solution $$\bfN(\lam;x,t) = (N_{11}(z;x,t),N_{12}(z;x,t)) $$ of RHP \ref{RHP2.row}, one recovers the solution $q(x,t)$ of \eqref{DNLS2} via the asymptotic formula
\begin{equation}
\label{DNLS2:q.recon.bis}
q(x,t) = \lim_{z \to \infty} 2iz N_{12}(z;x,t).
\end{equation}
As we will see, each $\lam_j = u_j + iv_j$ gives rise to a soliton of the form moving with velocity 
c=
$ -4u_j$.
Exact formulas for the soliton in terms of $\lam_j$ and $C_j$ may be found in \cite[Appendix B]{JLPS18b}.

We conclude this section by  introducing some additional notations. We denote by $U$ the subset of $H^{2,2}(\R)$ consisting of  functions $q$ for which $\balpha$ has no zeros on $\R$ and at most finitely many simple zeros in $\C^+$.
The set $U = \bigcup_{N=0}^\infty \, U_N$ where $U_N$ consists of functions $q$ for which $\balpha$ has exactly $N$ zeros in $\C^+$. If  $N \neq 0$, we denote by $\lam_1, \ldots, \lam_N$ the simple zeros of $\balpha$ in $\C^+$.

\subsection{Nonlinear steepest descent and long-time behavior}
Our  long-time behavior of solutions to DNLS is restricted to initial conditions $q_0$ in the subset of $H^{2,2}$
consisting of functions for which $\balpha$ has no zeros on $\R$ and at most finitely many simple zeros $\{\lambda_i\}_{1}^{N}$ in the upper-half complex plane  $\C^+$.  
Recall that these
zeros are referred to as resonances. They are responsible for the presence of separated solitons in the 
description of the long-time behavior of $q(x,t)$ and the soliton-resolution conjecture. 
Soliton resolution refers to the property that the solution decomposes into the sum of a finite number of separated solitons and a radiative part as 
$|t| \to \infty$.

The inverse scattering  method provides a full  description of  the asymptotic solution.  The limiting soliton parameters are slightly modulated, due to the soliton-soliton and soliton-radiation interactions.
 
Our analysis builds upon  the steepest descent method of Deift and Zhou \cite{DZ93, DZ03}, the later approach of  McLaughlin-Miller \cite{MM06} and  Dieng-McLaughlin \cite{DM08}, and the  work of Borghese, Jenkins and McLaughlin \cite{BJM18} on the focusing cubic NLS which shows how to treat a problem with discrete as well as continuous spectral data. (For more details on this approach, see the recent review article by Dieng, McLaughlin, and Miller \cite{DMM19}.)

Here, we  exclude initial conditions with spectral singularities. The latter may affect the long-time behavior of solutions in the same way that  resonances affect  the long-time behavior of solutions and this will be the object of a forthcoming work.

As we will explain, a direct consequence of the large-time asymptotics we obtain is the asymptotic stability (in $L^\infty$ norm) of $N$-soliton solutions.

\subsection{Statement of results}

We now state our results on global wellposedness, soliton resolution, and asymptotic stability of solitons for the DNLS equation.

\begin{theorem} (Global wellposedness)
\label{GWP}
Suppose that $q_0 \in H^{2,2}(\R)$. There exists a unique solution $q(x,t)$ of \eqref{DNLS2} with $q(x,t=0)=q_0$ and
$t \mapsto q(\dotarg,t) \in C([-T,T],H^{2,2}(\R))$ for every $T>0$. Moreover, the map $q_0 \mapsto q$ is Lipschitz continuous from $H^{2,2}(\R)$ to $C([-T,T], H^{2,2}(\R))$ for every $T>0$.
\end{theorem}

The main steps of the proof consist in a careful analysis of the direct and inverse map. The set of  scattering data $\calD$ depends on the properties of the
coefficient $\alpha(\lambda)$.  In the simple case of no resonances and no spectral singularities, the set  $\calD$ reduces to the 
reflection coefficient $\rho(\lambda)$. If $\alpha$ presents $N$  resonances $ \{\lambda_j\}_{j=1}^N$, 
then $\calD = \{ \rho, \{\lambda_j, C_j\}_{j=1}^N \}$, where the coefficients  $c_j$  defined above are  called the norming constants.
In the more complex case where spectral singularities are present, it is necessary to include in  $\calD$   
jump matrices for a Riemann-Hilbert problem defined on an extended contour. The entries of these matrices involve  Jost functions
associated to the initial potential $q_0$. This was a seminal idea introduced in a series of papers
by Zhou that we implemented in \cite{JLPS19} and allowed us to establish a global wellposedness result without restriction on initial conditions. A review of this approach is presented in Section \ref{Arbitrary}.   Pelinovsky et al \cite{PSS17} proved a global wellposedness result for initial
condtions supporting solitons but without spectral singularities.

In order to state our soliton resolution result, it is necessary to introduce some notations. 
To keep things as simple as possible we consider here the case where $\eps = -1$ and $t \to +\infty$. Let $\lambda_0 = -x/(4t)$.
Denote by $\calG$ the gauge transformation (cf.\ \eqref{DNLS:q.gauge} with $\eps=-1$)
\begin{equation}
\label{gauge-transform}
\left(\calG u \right)(x) = u(x) \exp\left( i \int_x^\infty |u(y)|^2 \, dy\right). 
\end{equation}
We will use $\calG$ to move back and forth between \eqref{DNLS} (the equation of interest) and \eqref{DNLS2} (the equation more convenient for inverse scattering).

Next, we define space-time `windows' for soliton resolution. Choose intervals $[v_1,v_2]$ of velocities and $[x_1,x_2]$ of initial positions.  We will  compute the asymptotic behavior of $q(x,t)$ in space-time regions of the form
\begin{equation}
\label{S.cone}
\calS(v_1,v_2,x_1,x_2)	=	
\left\{ (x,t) : x=x_0 + vt \text{ for } v \in [v_1,v_2], \, \, x_0 \in [x_1,x_2] \right\}.
\end{equation}

Recall that a soliton associated to eigenvalue $\lambda$ moves with velocity
($-4 \real  \lambda$).
Given an interval $I \subset \R$, we set
\begin{equation}
\label{LamI}
\Lambda(I) = \{ \lam \in \Lambda: \real (\lam) \in I \} = \Lambda^+(I) \cup \overline{\Lambda^+(I)}
\end{equation}
and
\begin{equation}
\label{NI}
N(I) = |\Lambda^+(I)|. 
\end{equation}
Solitons in $\Lambda([-v_2/4, -v_1/4])$ should be `visible' {within $\calS(v_1,v_2,x_1,x_2)$}, but remaining solitons will 
move either too slowly or too fast to be seen in the moving window. 
We also define
\begin{equation}
\begin{aligned}
\label{posint}
I^-	&=	\left\{ \lam \in \C:		
  \imag \lam = 0, \quad \real  \lam < \inf I \right\}, \\
I^+	&=	\left\{ \lam \in \C:		
  \imag \lam = 0, \quad  \real \lam > \sup I \right\}.
\end{aligned}  
\end{equation}

\begin{figure}[H]
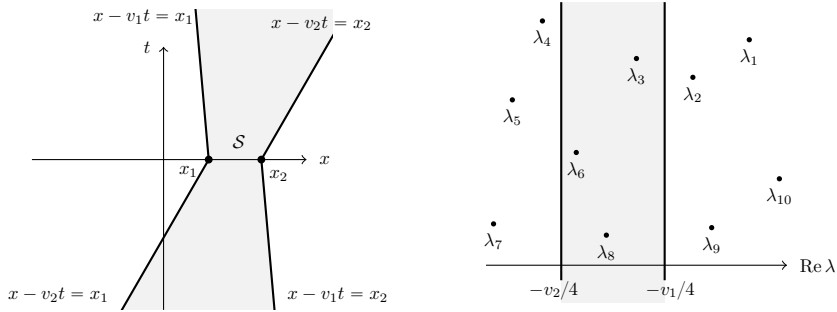

\solfig
\caption{
Given initial data $q_0(x)$ which generates scattering data 
$\left\{ \rho, \{ \lam_k, C_k \}_{k=1}^N \right\}$, then, asymptotically as $|t| \to \infty$ inside the space-time cone $\mathcal{S}(v_1,v_2,x_1,x_2)$ (shaded on left) the solution $u(x,t)$ of \eqref{DNLS} approaches an $N(I)$-soliton $\usol(x,t)$ corresponding to the discrete spectra in $\poles(I)$ 
(shaded region on right).  The connection coefficients $\widehat{C}_k$ for $\usol$ are modulated by the soliton-soliton and soliton-radiation interactions as described in Theorem~\ref{thm:soliton-resolution}.
}
\end{figure}

Finally, we define soliton solutions to \eqref{DNLS} via gauge-transformed soliton solutions of \eqref{DNLS2}. Let $u_0 \in H^{2,2}(\R)$ be given initial data for \eqref{DNLS} and let $q_0 = \calG (u_0)$. We suppose that $u_0$ is so chosen that $q_0 \in U$, i.e., $q_0$ supports at most finitely many solitons.  Denote by 
$\{ \lam_k, C_k \}$ the discrete scattering data of $q_0$, i.e., the zeros of $\alpha$ and the associated norming constants. For a given interval $I$, denote by  $\qsol(x,t;\mathcal{D}_I)$ the soliton solution of \eqref{DNLS2} with 
{modulating reflectionless scattering data}
$$ \calD_I = \left\{ \rho_I \equiv 0, \{(\lam_k, \widehat{C_k})\}_{\lam_k \in {\Lambda(I^+)}} \right\}$$
where
$$
\widehat{C_k} =  C_k \  
	\prod_{\mathclap{ \real \lam_j < -v_2/4} }
	\quad \ \left(	\frac{\lam_k - \lam_j}{\lam_k - \overline{\lam_j}} \right)^2 
	\exp \left(  
					\frac{i}{\pi} \int_{-\infty}^{\lambda_0} 
						\frac{\log\left( 1 + \lam |\rho(\lam)|^2 \right)}{\lam-\lam_k}  \, d\lam
			\right).
$$
and set
\begin{equation}
\label{u-soliton}
u_{\sol}(x,t) = \calG^{-1} \qsol(x,t). 
\end{equation}

We can now  give our result on long-time asymptotics. 
We also assume that $|\lambda_0| > M t^{-1/8}$ to keep the formulas uncluttered, and do not give an explicit formula for the dispersive ($O(t^{-1/2})$) term.  The full result for both signs of $t$ and arbitrary $\lambda_0$, together with explicit formulas for the dispersive term,  may be found in \cite[Theorem 1.6]{JLPS18b}.

\begin{theorem}
\label{thm:soliton-resolution}
Suppose that $u_0 \in H^{2.2}(\R^2)$,  and let $q_0$ be given by the gauge transformation \eqref{gauge-transform}. Suppose that $q_0 \in U$  and let $ \left\{ \rho, \{ (\lam_k,C_k) \}_{k=1}^N \right\}$ be the scattering data for $q_0$.  Fix $v_1,v_2,x_1,x_2$ with $v_1 < v_2$ and $x_1 < x_2$, and   
let $I=[-v_2/4,-v_1/4]$,   $\lambda_0 = -x/(4t)$.
Finally, fix $M > 0$ and assume that $|\lambda_0 | > Mt^{-1/8}$. 

The solution $u(x,t)$ of \eqref{DNLS} with $\eps=-1$ has the following asymptotics as $t \to \infty$ in the cone
$\calS(v_1,v_2,x_1,x_2)$:
$$ 
    u(x,t) =  
      u_{\sol}(x,t;\calD_I) e^{i\alpha_0(\lambda_0,+)}
      \left[ 1 +O(t^{-1/2})\right] 
$$
where
$$
	\alpha_0(\lambda_0,\pm)	= 	
		\pm \frac{1}{\pi} \int_{\mp \infty}^{\lambda_0} \frac{\log(1+ \lambda | \rho(\lam)|^2)}{\lam} \, d\lam
	  + 4 \sum_{\mathclap{\real  \lam_k \in I^\mp}} \arg \lam_k , 
$$	  
\end{theorem}

\begin{remark}
The phase $e^{i\alpha_0(\lam_0,+)}$ in the above formulas arises because of  the
mismatch in the phase of $u(x,t) = \calG^{-1}(q)$ and $\usol(x,t;\calD_I)$ as defined in \eqref{u-soliton}. The mismatch is  caused by the cumulative interaction of $\usol(y,t;\calD_I)$ with the radiation and soliton components of the full system which are traveling faster than our chosen reference frame $\calS$. Because the velocities are proportional to $- \real  \lam$ (recall that $v = - 4 \real  \lam$ for solitons), faster velocities correspond to the part of the spectrum $I^-$ which lies to the left of $\xi$ as $t\to \infty$.
\end{remark}	

Our analysis provides a proof of asymptotic stability of  $N$-soliton solutions. Recall the two-parameter family of $1$-soliton solutions of \eqref{DNLS} defined by \eqref{u 1sol}. A single point of discrete spectrum $\lambda = \nu + i \mu$ with norming constant $C$ corresponds to the specific one-soliton 
\[
	u_{\omega,c} (x-x_0 ,t ) e^{-i \varphi_0}
\] 
where 
$$ \omega = 4|\lambda|, \quad c = - 4 \nu, \quad x_0 = \frac{1}{4 \mu} \log \frac{ |\lam| |C|^2 }{4 \mu^2},$$ and 
$$\varphi_0 = \arg(\lambda) + \arg(C) + \pi/2.$$

%

\begin{theorem}
\label{asympt-stab}
Let $u_\mathrm{sol}(x,t)$ be an $N$ soliton solution of \eqref{DNLS} with
$u_{\mathrm{sol}}(x,0)  =    \calG^{-1} ( q_{\mathrm{sol}}(x,0) )   \in \calG^{-1}(U_N)$ and
with scattering data $\mathcal{D}^\mathrm{sol} = \{ 0, \{\lam_k^{\mathrm{sol}}, C_k^{\mathrm{sol}} \}_{k=1}^N \}$ such that $\real  \lam_k \neq \real  \lam_j$, $j\neq k$. There exist positive constants $\eta_0=\eta_0(q_{\mathrm{sol}})$, $T=T(q_{\mathrm{sol}})$, and $K= K(q_{\mathrm{sol}})$ such that any initial data $u_0 \in H^{2,2}(\R)$ with
\[
	\eta_1 := \| u_0 - u_{\mathrm{sol}}(\, \cdot\, , 0) \|_{H^{2,2}(\R)} 	\leq \eta_0
\]
also lies in $\mathcal{G}^{-1}(U_N)$ with scattering data 
$\mathcal{D} = \{ \rho, \{ \lam_k, C_k \}_{k=1}^N \}$ 
satisfying 
$$\|\rho\|_{H^{2,2}(\R)} +	\sum_{k=1}^N | \lam_k - \lam_k^\mathrm{sol}| + |C_k - C_k^\mathrm{sol}| \leq K \eta_1.
$$
Moreover, the solution of the Cauchy problem \eqref{DNLS} with initial data $u_0$  asymptotically separates into a sum of N 1-solitons 
\[
	\sup_{x \in \R} 
	\,\,
	\left| 
	\,
	u(x,t) - \sum_{k=1}^N  u_{\omega_k, c_k} (x -x_k^\pm,t ) 
	e^{ i( \alpha_0(\real  \lam_k, \pm) - \varphi_k^\pm) }
	\,
	\right|
	\leq K \eta_1 |t|^{-1/2}, 
	\quad   |t| > T
\]	
where the phase corrections $\alpha_0(\xi, \pm)$ are defined in Theorem~\ref{thm:soliton-resolution}, and
writing $\lam_k = \nu_k + i \mu_k$, the soliton parameters are given by $\omega_k = 4 |\lam_k|$, $c_k = -4\nu_k$, 
\begin{gather}
\label{x phases}
	x^\pm_k = 
		\frac{1}{4\mu_k} \log \left| \frac{\lam_k C_k^2}{4 \mu_k^2} \right|
		+ \frac{1}{2\mu_k} 
		\sum_{\substack{\lam_j \in \poles^+ \\ \mathclap{\pm (\nu_k - \nu_j) > 0}}} 
		  \log \left| \frac{ \lam_k - \lam_j}{\lam_k - \overline{ \lam_j} } \right|
			\mp \int_{\nu_k}^{\mp\infty}
			\frac{\kappa(s)}{(s-\nu_k)^2+\mu_k^2} ds \\
\label{alpha phases}
	\varphi^\pm_k = 
		\arg \lp i \lam_k C_k \rp 
		+ \sum_{\substack{\lam_j \in \poles^+ \\ \mathclap{\pm (\nu_k - \nu_j) > 0}}}
		    \arg \lp \frac{ \lam_k - \lam_j}{\lam_k - \overline{ \lam_j} } \rp
		  \pm 2 \int_{\nu_k}^{\mp \infty} 
		\frac{ (s - \nu_k) \kappa(s) }{(s-\nu_k)^2+\mu_k^2} ds.
		 \mod{2\pi}
\end{gather}
and $\kappa(s)$ is defined by \eqref{delta}. 
\end{theorem}

\section{The Direct and Inverse Maps}

 In this section we describe in greater detail the direct and inverse scattering maps. We begin with the direct map. Letting $\lambda=\zeta^2$, we get a new linear problem from (\ref{DNLS2:M})
\begin{subequations}
\label{n}
\begin{align}
\label{n.de}
\frac{d \bff{n}^\pm}{dx}	&=	-i\lambda \ad(\sigma_3) \bff{n}^\pm +\twomat{0}{q}{-\lambda \overline{q}}{0} \bff{n}^\pm + \bff{P} \bff{n}^\pm\\
\label{n.ac}
\lim_{x \rightarrow \pm \infty} \bff{n}^\pm(x,\lambda)	&=	\I
\end{align}
\end{subequations}
and the solutions are related by
\begin{equation}
\label{n.T}
\bff{n}^+(x,\lambda) =	\bff{n}^-(x,\lambda)
	e^{-i\lambda x \ad(\sigma_3)} 
	\twomat{\alpha(\lambda)}
				{\beta(\lambda)}
				{\lambda {\bbeta(\lambda)}}
				{{\balpha(\lambda)}},
\end{equation}
where $\bff{n}^\pm$ is a $2\times 2$ matrix:
\begin{equation}
\label{jost-n}
\bff{n}^\pm(x,\lambda)=\twomat{n^\pm_{11}(x,\lambda)}{n^\pm_{12}(x,\lambda)}{n^\pm_{21}(x,\lambda)}{n^\pm_{22}(x,\lambda)} .
\end{equation}

For the reader's convenience, we fix $\eps=-1$ throughout this section.
Zeros of $\alpha$ and $\balpha$ do occur for data of physical interest. Zeros of $\alpha$ and $\balpha$ on the real axis may occur and correspond to \emph{spectral singularities}. RHP \ref{RHP2.row} is no longer solvable since  the jump matrix now has singularities on  $\R$; moreover, any zeros of $\alpha$ and $\balpha$ in their domains of analyticity will make the Beals-Coifman solutions meromorphic rather than analytic. 
Thus, one of the key issues of inverse scattering transform is to properly treat eigenvalues and spectral singularities introduced by the zeros of $\alpha$ and $\balpha$. 

Here we will consider two distinct cases: the generic or $N$-soliton case,  when $\alpha$ and $\balpha$ have finitely many simple zeros in $\C \setminus \R$  (section \ref{N soliton}), and the general case where $\alpha$ and $\balpha$ may have countably many zeros that may accumulate on the real axis and additional zeros on the real axis (section \ref{Arbitrary}).  In the first case, we can establish global well-posedness,  soliton resolution, and asymptotic stability of solitons with precise asymptotics; in the second case, we can establish global well-posedness but cannot, at present, compute long-time asymptotics.

\subsection{The $N$-Soliton case}
\label{N soliton}
Setting $\lambda=\zeta^2$, we first note that  the evenness of $\ba$ and $a$ induces
$$\alpha(\lambda)=a(\zeta), \, ~ \balpha(\lambda)=\ba(\zeta).$$
$\balpha$ has analytic continuation into $\bbC^+$ while $\alpha$ has analytic continuation into $\bbC^-$
and they are given by the following Wronskians respectively:
$$\ba(\zeta)=\balpha(\zeta^2)=\begin{vmatrix}
    \psi_{11}^-& \psi_{12}^+ \\ 
    \psi_{21}^- & \psi_{22}^+ \\ 
    \end{vmatrix}, \, ~ a(\zeta)=\alpha(\zeta^2)=\begin{vmatrix}
    \psi_{11}^+& \psi_{12}^- \\ 
    \psi_{21}^+ & \psi_{22}^- \\ 
    \end{vmatrix}.$$
Now we illustrate the zeros of $\alpha$, $\balpha$ and $a$ and $\ba$ in their respective domain.
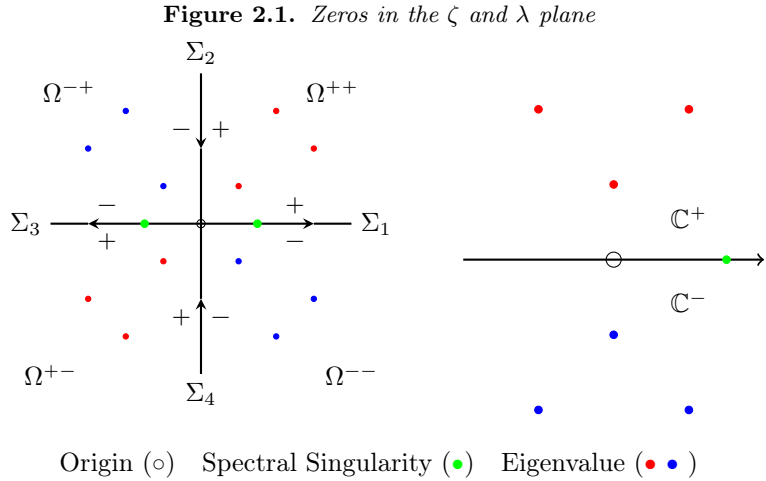
\begin{figure}[h!]
\caption{Zeros in the $\zeta$ and $\lambda$ plane}
\label{pic:zeta-lam}
\begin{tabular}{ccc}

\setlength{\unitlength}{5.0cm}
\begin{tikzpicture}[scale=0.5]
\draw[ thick,->,>=stealth] (0,0) -- (-3,0);
\draw[ thick] (-3,0) -- (-4,0);
\draw[thick,->,>=stealth] (0,0) -- (3,0);
\draw[ thick] (3,0) -- (4,0);
\draw[ thick,->,>=stealth] (0,4) -- (0,2);
\draw[ thick] (0,2) -- (0,0);
\draw[ thick,->,>=stealth] (0,-4) -- (0,-2);
\draw[ thick] 	(0,-2)--(0,0);
\node[above] at 		(-2.5,0) {$-$};
\node[below] at  		(-2.5,0) {$+$};
\node[above] at 		(2.5,0) {$+$};
\node[below] at 		(2.5,0) {$-$};
\node[left] at 			(0,2.5) {$-$};
\node[right] at 		(0,2.5) {$+$};
\node[left]  at 			(0,-2.5) {$+$};
\node[right] at 		(0,-2.5) {$-$};
\node[above] at (3.5, 3) {$\Omega^{++}$};
\node[above] at (-3.5, 3) {$\Omega^{-+}$};
\node[below] at (-4, -3.5) {$\Omega^{+-}$};
\node[below] at (4, -3.5) {$\Omega^{--}$};
\draw[black] (0,0) circle [radius=0.11];
\draw [green, fill=green] (1.5,0) circle [radius=0.10];
\draw [green, fill=green] (-1.5,0) circle [radius=0.10];
\draw [blue, fill=blue] (-3,2) circle [radius=0.07];
\draw [blue, fill=blue] (-2,3) circle [radius=0.07];
\draw [blue, fill=blue] (-1,1) circle [radius=0.07];
\draw [red, fill=red] (2,3) circle [radius=0.07];
\draw [red, fill=red] (3,2) circle [radius=0.07];
\draw [red, fill=red] (1,1) circle [radius=0.07];
\draw [red, fill=red] (-2,-3) circle [radius=0.07];
\draw [red, fill=red] (-3,-2) circle [radius=0.07];
\draw [red, fill=red] (-1,-1) circle [radius=0.07];
\draw [blue, fill=blue] (3,-2) circle [radius=0.07];
\draw [blue, fill=blue] (2,-3) circle [radius=0.07];
\draw [blue, fill=blue] (1,-1) circle [radius=0.07];
\node[right] at (4 , 0) {$\Sigma_1$};
\node[left] at (-4 , 0) {$\Sigma_3$};
\node[above] at (0 , 4) {$\Sigma_2$};
\node[below] at (0, -4) {$\Sigma_4$};

\end{tikzpicture}


&\hspace{1cm}&

\setlength{\unitlength}{5.0cm}
\begin{tikzpicture}
\draw [red, fill=red] (0,1) circle [radius=0.05];

\draw [red, fill=red] (-1,2) circle [radius=0.05];

\draw [red, fill=red] (1,2) circle [radius=0.05];

\draw [black]  (0,0) circle [radius=0.1];
\node [above] at (1,0.3) {$\mathbb{C^+}$};
\draw [->][thick] (-2,0) -- (2,0);
\node [below] at (1,-0.3) {$\mathbb{C^-}$};
\draw [blue, fill=blue] (0,-1) circle [radius=0.05];

\draw [green, fill=green] (1.5,0) circle [radius=0.05];
\draw [blue, fill=blue] (1,-2) circle [radius=0.05];

\draw [blue, fill=blue] (-1,-2) circle [radius=0.05];

\end{tikzpicture}

 \\[0.2cm]

\end{tabular}

\vskip 0.1cm

\begin{tabular}{ccc}
Origin ({\color{black} $\circ$}) &
Spectral Singularity ({\color{green} $\bullet$})	&	
Eigenvalue ({\color{red} $\bullet$} {\color{blue} $\bullet$} ) 
\end{tabular}

\label{fig:spectra}
\end{figure}

Suppose that $\balpha$ has $\lbrace \lam_k \rbrace_{k=1}^N$ simple zeros in $\bbC^+$. Then
\begin{align}
\label{B_0}
   \begin{bmatrix}
           n^-_{11}(x,\lambda_k) \\
           n^-_{21}(x,\lambda_k) \\
        \end{bmatrix}=B_k\lambda_k \begin{bmatrix}
           n^+_{12}(x,\lambda_k) \\
           n^+_{22}(x,\lambda_k) \\
        \end{bmatrix}e^{2ix{\lambda_k}}
  \end{align}
and
\begin{align}
\label{B_k*}
   \begin{bmatrix}
           n^-_{12}(x,\overline{\lambda}_k) \\
           n^-_{22}(x,\overline{\lambda}_k) \\
        \end{bmatrix}=-\overline{B}_k\begin{bmatrix}
           n^+_{11}(x,\overline{\lambda}_k) \\
           n^+_{21}(x,\overline{\lambda}_k) \\
        \end{bmatrix}e^{-2ix{\overline{\lambda}_k}}.
  \end{align}
 We also set \begin{equation}
 \label{norming-lam}
 C_k=\dfrac{B_k}{ \balpha'(\lambda_k)}.
 \end{equation}
 and define the discrete scattering data
 \begin{equation}
 \label{discrete data}
 \lbrace\lam_k, C_k\rbrace_{k=1}^N.
 \end{equation}
 
 Scattering data evolve linearly in time: \begin{subequations}
 \begin{align}
 \label{time-rho}
 \dot\rho(\lam,t)&= -4i \lam^2 \rho(\lam,t)\\
 \label{time-lambda}
\dot\lam_k &=0,\quad k=1,...,N \\
\label{time-C}
\dot C_k &=-4i\lam_k^2C_k,
 \end{align}
\end{subequations}
Thus
\begin{align}
 \begin{cases}
\rho(\lam,t)&= e^{-4i \lam^2t} \rho(\lam)\\
\lam_k (t)&=\lam_k,\quad k=1,...,N \\
 C_k(t) &=e^{-4i\lam_k^2 t}C_k ,
\end{cases}
\end{align}
where $\rho$, $\lambda_k$ and $C_k$ are associated to the initial data.
The proof of the following theorem can be found in \cite[Section 3]{JLPS18a}.
\begin{proposition}

\label{thm:R}
There is a spectrally determined open and dense subset $U=\bigcup_{n=0}^\infty U_n$ of $H^{2,2}(\bbR)$ containing a neighborhood of $0$ so that for $N=0,1,...$ the direct scattering map $\mathcal{R}$ 
\begin{align*}
\mathcal{R}: U_N &	\longrightarrow	H^{2,2}(\bbR) \times (\bbC_\times\times \bbC^+)^N \\
q							&\mapsto				(\rho, \lbrace  C_i, \lambda_i \rbrace_{i=1}^N)
\end{align*}
is a  Lipschitz continuous map from bounded subsets of $U_N$ into bounded subsets of $V_N=H^{2,2}(\R) \times (\bbC_\times \times \bbC^+)^N$.
\end{proposition}
We list the main steps of proof here.
\begin{itemize}
\item[1.] Write $\alpha$ and $\beta$ in terms of the Jost solutions
\begin{align}
\label{alpha}
\alpha(\lam)	&=	n_{11}^+(0,\lam) \overline{n_{11}^-(0,\lam)} + \lam^{-1}\overline{ n_{21}^-(0,\lam) } n_{21}^+(0,\lam) ,\\[10pt]
\label{beta}
\beta(\lam)		&=	\frac{1}{\lam}\left(- \overline{ n_{11}^-(0,\lam)} \overline{n_{21}^+(0,\lam)} 
+\overline{n_{11}^+(0,\lam)} \overline{n_{21}^-(0,\lam)}\right)
\end{align}
\item[2.]  Show that $\alpha$ and $\beta$ belong to appropriate  function spaces by studying the following two Volterra integral equations:
\begin{equation}
\nonumber
n_{11}^\pm(x,\lam)	
	=	1	-	 \int_x^{\pm\infty} 
								q(y) n_{21}^+(y,\lambda)  
						-	
							\frac{i}{2} |q(y)|^2 n_{11}^\pm(y,\lam) 
					\, dy 
\end{equation}	
\begin{equation}	
\nonumber			
n_{21}^+(x,\lam)=	\int_x^{\pm\infty}
				e^{2 i\lam(x-y)}   
				\left(		 \lam \overline{q(y)} n_{11}^\pm(y,\lam) 
				+  \frac{i}{2} 
				   |q(y)|^2 n_{21}^\pm(y,\lam)\right)dy.	
\end{equation}
 For more details, see \cite[Chapter 3.1]{Liu17}.
\item[3.] $\rho=\beta/\alpha\in H^{2,2}$ will follow from an application of the quotient rule.
\item[4.]Using the relations given in \eqref{B_0}-\eqref{norming-lam} and \eqref{alpha} to show that $\lbrace \lambda_k, C_k  \rbrace_{k=1}^N$ depends continuously on the initial data $q$. Detailed proofs are given in  \cite[Chapter 3.2-3.4]{Liu17}
\end{itemize}

The next step is to construct the potential $q$ from the above scattering data. More precisely,  $q$ is  reconstructed through  the solution of  RHP \ref{RHP2.row} which is in fact, equivalent to a Riemann-Hilbert problem  (RHP  \ref{RHP-2}  below) with no discrete data but having an augmented contour 
\begin{align}\label{Gamma} 
	\Gamma= \bbR \cup \{ \Gamma_j \}_{j=1}^{N} \cup  \{ \Gamma_{j}^* \}_{j=1}^{N}  
\end{align}
where each $\Gamma_j$ (resp. $\Gamma^*_j$)  is a simple closed curve in $\C \setminus \bbR$ surrounding $\lambda_j$ (resp. $\lambda_j^*$), as shown in Fig. \ref{figure-2}. 
The new curves are given an orientation consistent with the orientation of the original contour $\bbR$, so $\Gamma$ also divides $\C \setminus \bbR$ into 
two disjoint sets. 

\begin{figure}[H]
\caption{The Augmented Contour $\Gamma$}
\begin{center}
\label{figure-2}
\begin{tikzpicture}[scale=0.8]
\draw [red, fill=red] (0,1.5) circle [radius=0.05];
\draw[thick,->,>=stealth] (0.3,1.5) arc(360:0:0.3);
\draw [red, fill=red] (-1,3) circle [radius=0.05];
\draw[thick,->,>=stealth] (-0.7,3) arc(360:0:0.3);

\draw [red, fill=red] (1,3) circle [radius=0.05];
\draw[thick,->,>=stealth] (1.5,3) arc(360:0:0.5);

\node [below] at (1.45, 2.55) { $\Gamma_i$};
\node[above] at(1.15,2.95){\footnotesize{$-$}};
\node[above] at(1.55,3.25){\footnotesize{$+$}};

\draw  (0,0) circle [radius=0.1];
\node [above] at (1.4,0.3) {$\mathbb{C^+}$};
\draw [very thick,->,>=stealth] (-4,0) -- (4,0);
\node [below] at (1.4,-0.3) {$\mathbb{C^-}$};

\draw [blue, fill=blue] (-1,-3) circle [radius=0.05];
\draw[thick,->,>=stealth] (-0.7,-3) arc(0:360:0.3);

\draw [blue, fill=blue] (0,-1.5) circle [radius=0.05];
\draw[thick,->,>=stealth] (0.3,-1.5) arc(0:360:0.3);

\draw [blue, fill=blue] (1,-3) circle [radius=0.05];
\draw[thick,->,>=stealth] (1.5,-3) arc(0:360:0.5);
\node [above] at (1.45, -2.65) { $\Gamma_i^*$};
\node[below] at	(1.15,-2.95){\footnotesize{$+$}};
\node[below] at	(1.55,-3.25){\footnotesize{$-$}};

\end{tikzpicture}
\end{center}
\end{figure}

\begin{RHP}
\label{RHP-2}
Fix $x \in \bbR$ and   let $( \rho,  \{ C_i, \lambda_i  \}_{i=1}^{N} )\subset H^{2,2}(\bbR)\times (\bbC_\times\times \bbC^+)^N $. 
Find a vector-valued  function  $\bfN(\dotarg;x)$ with the following properties:
\begin{itemize}
\item[(i)](Analyticity) $\bfN(z; x)$ is a  row vector-valued  analytic function of $z$ for $z\in \mathbb{C}\setminus\Lambda $  where $\Gamma$, defined by \eqref{Gamma}, is depicted in Figure~\ref{figure-2}.

\smallskip

\item[(ii)] (Normalization) ${\bfN}(z;x)=(1 ,0)+\mathcal{O}(z^{-1})$ as $z \rightarrow\infty$.

\smallskip

\item[(iii)] (Jump condition) For each $\lambda\in\Gamma $, $\bfN$ has continuous 
boundary values
$\bfN_{\pm}(\lam;x)$ as $z \to \lambda$ from $\bbC^\pm$. 
Moreover, the jump relation 
$$\bfN_+(\lambda;x)=\bfN_-(\lambda;x) \bff{J}_{x}(\lambda)$$ 
holds, where for $\lambda\in\mathbb{R}$
$$
\bff{J}_{x}(\lambda)=
		e^{-i\lambda x\, \ad \sigma_3} 
		\Twomat	{1+\lambda|\rho(\lambda)|^2}{\rho(\lambda)}
						{\lambda\overline{\rho(\lambda)}}{1}
$$

\smallskip
\item[(iv)]  (Residue condition)
For each $\lam \in \Gamma_i\cup\Gamma_i^* $ 
$$
\bff{J}_x(\lambda) = 	\begin{cases}
						\twomat{1}{0}{\dfrac{C_i \,\lambda_i e^{2i\lambda x}}{\lambda-\lambda_i}}{1}	&	\lambda \in \Gamma_i, \\
						\\
						\twomat{1}{\dfrac{\overline{C_i} \, e^{-2ix \lambda} }{\lambda-\overline{\lambda_i}}}{0}{1}
							&	\lambda \in \Gamma_i^*
					\end{cases}
$$
\end{itemize}
\end{RHP}
One recovers $q$ from the relation
\begin{align*}
 q(x) = 2i\lim_{z \rightarrow\infty} z\, {\mathbf{N}}_{12}(z;x)
 \end{align*}
where the limit is 
 uniform as $|z| \to \infty$  in proper subsectors of $\C \setminus \R$.

Beals and Coifman showed that Riemann-Hilbert problems such as Problem \ref{RHP-2} can be reduced to a boundary integral equation on the jump contour  known as the Beals-Coifman integral equation, much as a boundary value problem for an elliptic equation can be reduced to a boundary integral equation. 

In our case, the Beals-Coifman integral equation takes the following form. First, observe that the 
jump matrix $$\bff{J}_x(\lambda) \coloneqq e^{i\lam x \ad \sigma_3} \bff{J}(\lam)$$ admits the  factorization
$$
\bff{J}_x(\lam) = (\I-\bff{W}_x^-)^{-1} (\I + \bff{W}_x^+)
$$
where
$$
\left( \bff{W}_x^-, \bff{W}_x^+ \right) =
\begin{cases}
\left(
	\twomat{0}{0}{0 \vphantom{\frac{\lam_i}{\lam-\lam_i}}}{0} 
	\, , \,
	\twomat{0}{0}{\dfrac{\lambda_i C_i e^{2ix\lambda} }{\lambda- \lambda_i}}{0} 
\right)	&	\lambda \in \Gamma_i,\\
\\
\left(
	\twomat{0}{-\dfrac{\overline{C_i} e^{-2ix\lam}}{\lam-\overline{\lam_i}}}{0}{0}
	\, , \,
	\twomat{0}{0 \vphantom{\frac{1}{\lam-\lam_i}}}{0}{0}
\right)	&	\lambda \in \Gamma_i^*, \\
\\
\left(
	\twomat{0}{-\rho(\lam) e^{-2ix\lam}}{0}{0}
	\, , \,
	\twomat{0}{0}{\lam \overline{\rho(\lam)} e^{2ix\lam}}{0}
\right), &	\lam \in \R.
\end{cases}
$$
Next, denote by $C^+_\Gamma$, $C^{-}_\Gamma$ the Cauchy projection operators for the contour $\Gamma$. Finally, define a new unknown row vector-valued function $\bff{\nu}$ on $\Gamma$ by
$$ \bff{\nu}(\lam;x) =  \bfN_+(\lam;x) (\I + \bff{W}_x^+)^{-1} = \bfN_-(\lam;x)(\I-\bff{W}_x^-)^{-1}. $$ 
Note that 
$$ \bfN_+(\lam;x) - \bfN_-(\lam;x)  = \bff{\nu}(\lam;x) \left(\bff{W}_x^-(\lam) + \bff{W}_x^+(\lam)\right) $$
gives the jump of $\bfN$ across $\Gamma$, so that, given $\bff{\nu}(\lam;x)$, one can recover
$\bfN(z;x)$ via the Cauchy integral
$$\bfN(z;x) =(1,\, 0 )+ \frac{1}{2\pi i} \int_{\Gamma} \frac{ \bff{\nu}(z;x) \left(\bff{W}_x^+(s)+\bff{W}_x^-(s) \right)}{s-z} \, ds.$$
A standard argument using this Cauchy integral representation leads to the 
\emph{Beals-Coifman integral equation} for the function $\bff{\nu}$:
\begin{align}
\label{SIE-2}
\bff{\nu}(x, \lambda)&=(1,\, 0)+\mathcal{C}_w  \bff{\nu}(\lambda;x)\\
\intertext{where $\calC_w$ is the \emph{Beals-Coifman integral operator} acting on a row vector-valued function $\bff{f}$ on $\Gamma$ by}
\nonumber
   \left(\mathcal{C}_w \bff{f} \right)(\lam)
   	&:=C^+_\Gamma(\bff{f} \bff{W}_x^-)(\lambda)+C^-_\Gamma (\bff{f} \bff{W}_{x}^+)(\lambda).
\end{align}

The solvability of RHP \ref{RHP-2} follows from the analysis and solvability of a similar RHP \cite[RHP 5.1.4]{Liu17} in the original $\zeta$ variable  \cite[Lemma 5.2.3]{Liu17}.
We finally obtain the reconstructed potential
\begin{align*}
\label{q.recon.nu}
q(x) &= 
		-\frac{1}{\pi} \int_{-\infty}^\infty {{\bff{\nu}}}_{11}(s;x) \rho(s) e^{-2isx} ds 
		-\sum_{k=1}^n  2i  {\bff{\nu}}_{11}(\overline{\lam_k};x) \overline{C_k} e^{-2i\overline{\lam_k}x}\\
		\end{align*}

\subsection{Arbitrary Spectral Singularities}
\label{Arbitrary}
In this section we allow $\alpha(z)$ and $\balpha(z)$ to have (possibly infinitely many) zeros  of arbitrary order in $\bbC^\pm\cup \bbR$. This leads to two consequences:
\begin{enumerate}
\item The reflection coefficients $\rho$ are not defined on all of $\bbR$.
\smallskip
\item The Beals-Coifman solution $\bff{N}(z;x)$ may not have  a continuous limit as $z$ approaches $\bbR$ from $\bbC^\pm$.
\end{enumerate}
In a series of papers \cite{Zhou89-1, Zhou89-2, Zhou98},  Zhou developed new tools to construct  direct and inverse scattering maps 
in the presence of arbitrary spectral singularities.
The key idea is to make use of the following two 
observations about
$\balpha(z)$:
\begin{enumerate}
\item For $\left\Vert    q    \right\Vert_{L^1} \ll 1$, $\balpha$  has no zeros in $\bbC^+\cup \bbR$. This follows from 
analytic Fredholm theory.
\smallskip
\item The asymptotic relation $\lim_{z \to \infty} \balpha(z) =1$  holds. Thus all the zeros of $\alpha$ and $\balpha$ are located in a bounded region of $\bbC$. This follows from normalization and the Riemann-Lebesgue lemma.
\end{enumerate}
\begin{figure}[H]
\caption{The Augmented Contour $\Gamma$}
\bigskip
\begin{tikzpicture}[scale=0.9]
\draw  (0,0) circle [radius=0.1];
\draw	 [blue, fill=blue] (3,0) circle [radius=0.1];
\draw [blue, fill=blue] (-3,0) circle [radius=0.1];
\draw	 [thick,->,>=stealth] (-3,0) arc(180:90:3);
\draw [thick,->,>=stealth] (-3,0) arc(180:270:3);
\draw[thick]	(0,3) arc(90:0:3);
\draw[thick]	(0,-3) arc(270:360:3);
\draw[thick] (0,0) circle [radius=3];
\draw [thick,->,>=stealth] (-5,0) -- (-4,0);
\draw [-][thick] (-4,0) -- (-3,0);
\draw [->,>=stealth][thick] (3,0) -- (4,0);
\draw[-][thick] (4,0) -- (5,0);
\draw [<-,>=stealth][thick] (-2,0) -- (3,0);
\draw[-][thick] (-2,0) -- (-3,0);
\node [below] at (1,-0.3) {$+$};
\node [above] at (1,0.3) {$-$};
\node [below] at (4,-0.3) {$-$};
\node [above] at (4,0.3) {$+$};
\node [below] at (-4,-0.3) {$-$};
\node [above] at (-4,0.3) {$+$};
\node [below] at (0,-0.5) {$\Omega_4$};
\node [above] at (0,0.5) {$\Omega_3$};
\node [below] at (4,-1) {$\Omega_2$};
\node [above] at (4, 1) {$\Omega_1$};
\node [above] at (-2.8, 2) {$\Sigma_\infty$};
\node [right] at (5,0) {$\mathbb{R}$};
\node[below] at (2.6,0){$S_\infty$};
\node[below] at (-2.4,0){$-S_{\infty}$};
\node [above] at (2.8, 2) {$\Sigma_\infty^+$};
\node [below] at (2.8, -2) {$\Sigma_\infty^-$};
\end{tikzpicture}
\label{pic:sing}
\end{figure}
\begin{remark}
The regions $$\Omega_+ = \Omega_1 \cup \Omega_4, \quad \Omega_- = \Omega_2 \cup \Omega_3$$ lie, respectively, to the left and right of $\Gamma$. We set 
$$\Gamma_+=\partial\Omega_1 \cup \partial \Omega_4. \quad \Gamma_-=\partial\Omega_2 \cup \partial \Omega_3.$$ 
\end{remark}
\medskip

To make use of the  two observations
above, we first augment the real line with a circle $\Sigma_\infty$ that contains all the zeros of $\alpha$ and $\balpha$ and denote by $ \Gamma=\bbR\cup \Sigma_\infty$ the resulting augmented contour.
 The oriented contour $\Gamma$
separates the complex plane into $\Omega_1-\Omega_4$ as described in Figure \ref{pic:sing}. 


The construction of the scattering data  for the augmented contour has to be done in such a way that the  resulting inverse map is well-defined.
Let $x_0\in \bbR$ be such that the cut-off potential   $$q_{x_0}=q\chi_{ (x_0, \infty )  }$$ satisfies $ \left\Vert q_{x_0}\right\Vert_{L^1}  \ll 1$. We construct a new Beals-Coifman function inside the circle which is also normalized at $x\to +\infty$ to replace $\bfN(z;x)$ inside the circle $\Sigma_\infty$ and formulate a new Riemann-Hilbert problem along the augmented contour $\Gamma$. 
We can still construct $\bfN(z;x)$ and $\rho$ outside the circle $\Sigma_\infty$ as before,  given by Problem \ref{RHP-2}. We then construct a Beals-Coifman solution $\bfN^{(0)}$ normalized as $x\to \infty$ associated to the  potential $q_{x_0}$. $\bfN^{(0)}$ and $\bfN(z;x)$ are constructed simultaneously.  We define the scattering data across different parts of the contour $\R\cup \Sigma_\infty$ as illustrated in Figure \ref{fig-scattering}.

\begin{figure}[H]
\caption{Scattering data for $q$}

\bigskip

\begin{tikzpicture}[scale=0.9]
\draw[black, fill=black] 		(3,0) 		circle [radius=0.05];
\draw [black, fill=black] 		(-3,0) 	circle [radius=0.05];
\draw[thick,->-]		(-3,0) arc(180:0:3);
\draw[thick,->-]		(-3,0)	 arc(180:360:3);
\draw [thick,->-] 		(-6.5,0) -- (-3,0);
\draw [thick,->-]		(3,0)	--	(-3,0);
\draw [thick,->-] 		(3,0) -- 	(6.5,0);
\node [below] at (0,-2) 		{$\Omega^+$};
\node [above] at (0,2) 		{$\Omega^-$};
\node [below] at (4,-3) 		{$\Omega^-$};
\node [above] at (4, 3) 		{$\Omega^+$};
\node[below] at (2.5,0)		{$S_\infty$};
\node[below] at (-2.4,0)		{$-S_{\infty}$};
\node[above] at 	(-4.75,0)		{\footnotesize{$\bff{J}=\twomat{1+\lam|\rho|^2}{\rho}{\lam \overline{\rho}}{1}$}};
\node[above] at 	(4.75,0)		{\footnotesize{$\bff{J}=\twomat{1+\lam|\rho|^2}{\rho}{\lam \overline{\rho}}{1}$}};
\node[above]	at	(0,0)			{\footnotesize{$\bff{J}=\twomat{1}{-\rho_0}{-\lam \overline{\rho_0}}{1+\lam|\rho_0|^2}$}};
\node[above]	at	(0,3.2)			{\footnotesize{$\bff{J}=\twomat{1}{0}{e^{-2ix_0 \lam}\frac{n_{21}^-(x_0,\lam)}{\balpha \balpha_0}}{1}$}};
\node[below] 	at	(0,-3.2)			{\footnotesize{$\bff{J}=\twomat{1}{-e^{2ix_0\lam}\frac{n_{12}^-(x_0,\lam)}{\alpha \alpha_0}}{0}{1}$}};
\end{tikzpicture}
\label{fig-scattering}
\end{figure}
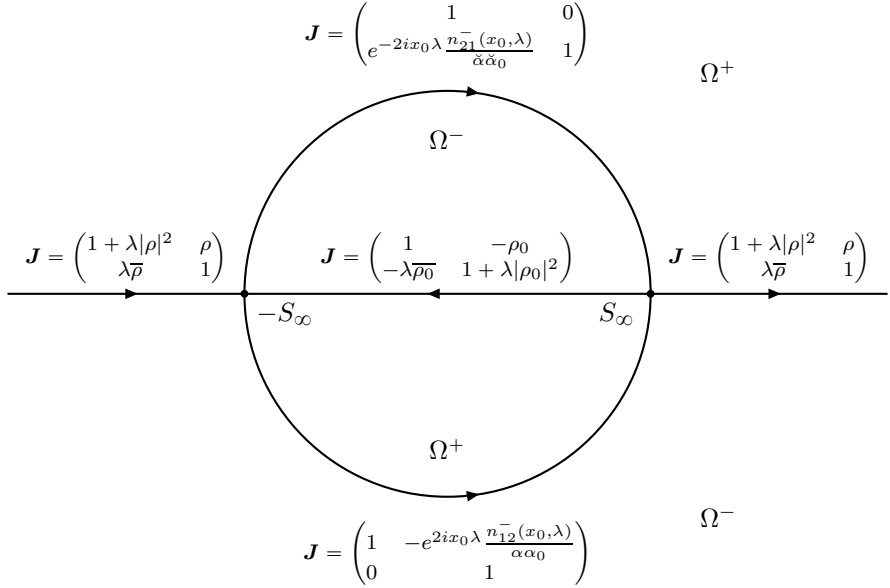
\begin{remark}
In Figure \ref{fig-scattering}, $n^-_{21}(x_0, \lambda)$ and $n^-_{12}(x_0, \lambda)$ refer to the Jost functions in \eqref{jost-n} evaluated at $x=x_0$. The quantitites $\rho_0$, $\alpha_0$ and $\breve{\alpha}_0$ are  scattering coefficients associated to $q_{x_0}$.
\end{remark}

We define the Sobolev spaces $H^k(\Gamma_\pm)$  needed for the study of arbitrary spectral singularities. These spaces were  introduced by Zhou \cite{Zhou89-1}.
If $\Gamma = \Gamma_1 \cup \ldots \cup \Gamma_n$ and the $\Gamma_i$ are either half-lines, line segments, or arcs, the space $H^k(\Gamma)$ consists of  functions $f$ on $\Gamma$ with the property that $\left. f \right|{\Gamma_i} \in H^k(\Gamma_i)$.  Limits of $f^{(j)}$ at the endpoints of  $\Gamma_i$ are well-defined for $0 \leq j \leq k-1$. The space $H^k(\Gamma_+)$ (resp.\ $H^k(\Gamma_-)$) consists of  the  functions of $H^k(\Gamma_i)$ which are continuous together with their derivatives up to order $k-1$ along the boundary of the positive (negative) components shown in Figure \ref{fig-Gamma+-comp} (resp.\ Figure \ref{fig-Gamma--comp}). 

\begin{center}
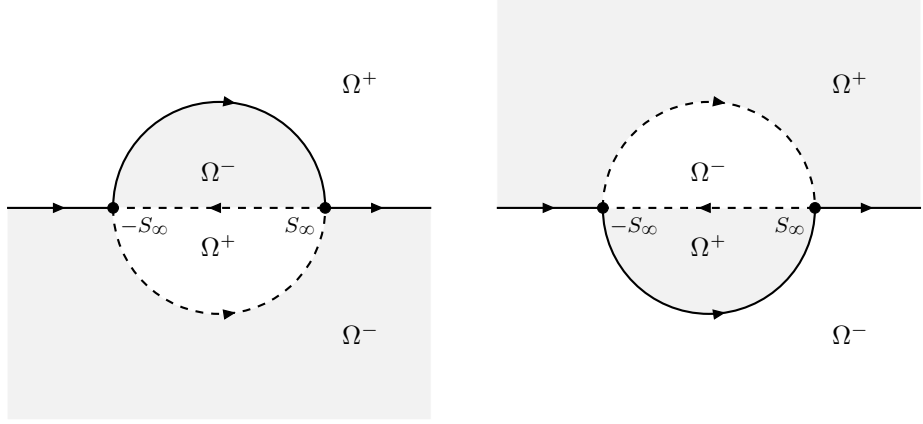
\begin{figure}[H]
\begin{subfigure}{0.45\textwidth}
\caption{Boundary components of $\Omega^+$}
\bigskip
\begin{tikzpicture}[scale=0.47]
\draw[white]						(-6,0)		--	(-6,6)		--	(6,6)	--	(6,0);
\draw[white,fill=gray!10]		(-6,0)		--	(-3,0)		-- 	(-3,0) arc(-180:0:3)	--	(6,0)	--
										(6,-6)		--	(-6,-6)	--	(-6,0);
\draw[white,fill=gray!10]		(-3,0)	 arc(180:0:3)	--	(-3,0);
\draw[black, fill=black] 		(3,0) 		circle [radius=0.15];
\draw [black, fill=black] 		(-3,0) 	circle [radius=0.15];
\draw[thick,->-]		(-3,0) arc(180:0:3);
\draw[thick,dashed,->-]		(-3,0)	 arc(180:360:3);
\draw [thick,->-] 					(-6,0) -- (-3,0);
\draw [thick,->-,dashed]		(3,0)	--	(-3,0);
\draw [thick,->-] 					(3,0) -- 	(6,0);
\node [below] at (0,-0.5) 		{$\Omega^+$};
\node [above] at (0,0.5) 		{$\Omega^-$};
\node [below] at (4,-3) 		{$\Omega^-$};
\node [above] at (4, 3) 		{$\Omega^+$};
\node[below] at (2.3,0)		{\footnotesize{$S_\infty$}};
\node[below] at (-2.1,0)		{\footnotesize{$-S_{\infty}$}};
\end{tikzpicture}
\label{fig-Gamma+-comp}
\end{subfigure}
\qquad
\begin{subfigure}{0.45\textwidth}
\caption{Boundary components of $\Omega^-$}
\bigskip
\begin{tikzpicture}[scale=0.47]
\draw[white]						(-6,0)		--	(-6,-6)		--	(-6,6)		--	(-6,0);
\draw[white,fill=gray!10]		(-6,0)		--	(-3,0)		-- 	(-3,0) arc(180:0:3)	--	(6,0)	--
										(6,6)		--	(-6,6)		--	(-6,0);
\draw[white,fill=gray!10]		(-3,0)	 arc(-180:0:3)	--	(-3,0);
\draw[black, fill=black] 		(3,0) 		circle [radius=0.15];
\draw [black, fill=black] 		(-3,0) 	circle [radius=0.15];
\draw[thick,->-,dashed]		(-3,0) arc(180:0:3);
\draw[thick,->-]					(-3,0)	 arc(180:360:3);
\draw [thick,->-] 					(-6,0) -- (-3,0);
\draw [thick,->-,dashed]		(3,0)	--	(-3,0);
\draw [thick,->-] 					(3,0) -- 	(6,0);
\node [below] at (0,-0.5) 		{$\Omega^+$};
\node [above] at (0,0.5) 		{$\Omega^-$};
\node [below] at (4,-3) 		{$\Omega^-$};
\node [above] at (4, 3) 		{$\Omega^+$};
\node[below] at (2.3,0)		{\footnotesize{$S_\infty$}};
\node[below] at (-2.1,0)		{\footnotesize{$-S_{\infty}$}};
\end{tikzpicture}
\label{fig-Gamma--comp}
\end{subfigure}
\caption{Boundary Components of $\Omega^\pm$}
\end{figure}
\end{center}

Let's now consider the following RHP:

\begin{RHP}
\label{RHP.lambda}
Fix $x \in \bbR$. Find a row vector-valued  function 
$\bfN(\dotarg; x)$ on $\C \setminus \Gamma$ with the following properties:
\begin{enumerate}
\item[(i)](Analyticity) $\bfN(z;x)$ is an analytic function of $z$ for $z\in \mathbb{C}\setminus\Gamma $,
\item[(ii)] (Normalization) ${\bfN}(z;x)=(1 ,0)+\bigO{z^{-1}}$ as $z \rightarrow\infty$,
and
\item[(iii)] (Jump condition) For each $\lambda\in\Gamma $, $\bfN$ has continuous 
boundary values
$\bfN_{\pm}(\lambda)$ as $z \to \lambda$ from $\Omega_\pm$. 
Moreover, the jump relation 
$$\bfN_+(\lambda;x)=\bfN_-(\lambda;x) \bff{J}_x(\lambda)$$ 
holds, where 
\begin{align*}
\bff{J}_x(\lambda)
	&=	
		\begin{cases}
				\ttwomat{ 1+\lambda|\rho(\lambda)|^2}{\rho(\lambda)e^{-2i\lam x}}{\lambda\overline{\rho(\lambda)}e^{2i\lam x}}{1},
				&	\lambda\in \R_\infty
				\\	
				\\
				\ttwomat{1} {-\rho_0(\lambda)e^{-2i\lam x}}{-\lambda\overline{\rho_0(\lambda)}e^{2i\lam x}}  { 1+\lambda|\rho_0(\lambda)|^2},
				&	\lam \in(-S_\infty, S_\infty)
				\\ \\
				\ttwomat{1}{0}{e^{-2ix_0\lambda}  \dfrac{n^{-}_{21} (x_0, \lambda) }{\balpha(\lambda) \balpha_0(\lambda)  }}{1}
				& \lambda \in \Sigma_\infty^+, \\
				\\
				\ttwomat{1}{-e^{2ix_0\lambda}\dfrac{n^{-}_{12} (x_0, \lambda) }{\alpha(\lambda) \alpha_0(\lambda)  }}{0}{1} 
				&	\lambda \in  \Sigma_\infty^-.
		\end{cases}
\end{align*}

\medskip

The  solvability of the RHP \ref{RHP.lambda} follows two main ingredients:

\medskip

\begin{enumerate}
\item[1.] The analysis of a Beals-Coifman 
integral equation similar to Equation \eqref{SIE-2}, and
\smallskip
\item[2.] A key triangular factorization of the jump matrices $\bff{J}(\lambda)$ as characterized in the following theorem.
\end{enumerate}

\begin{theorem} \cite[Theorem 2.11]{JLPS19}
\label{thm:scattering data}
 The jump matrix $\bff{J}(\lambda)$ represented  in  Figure \ref{fig-scattering}  along the different sections of the contour, admits a triangular factorization $$\bff{J}(\lambda)=\bff{J}_-^{-1}(\lambda) \bff{J}_+(\lambda)$$ where: 
\begin{enumerate}
\item[(i)]
$\bff{J}_{-}(\lambda)-\I\in H^{2,2}(\partial \Omega_2) $, $\bff{J}_{-}(\lambda)-\I\in H^{2}(\partial \Omega_3) $,
$\bff{J}_{+}(\lambda)-\I\in H^{2}(\partial \Omega_4)$ and $\bff{J}_+(\lambda)-\I\in H^{1,1}(\partial\Omega_1)$, and
\smallskip
\item[(ii)] $\bff{J}_+\restriction_{\partial\Omega_1}-\I$ and $\bff{J}_-\restriction_{\partial\Omega_3}-\I$ are strictly lower triangular while $\bff{J}_-\restriction_{\partial\Omega_2}-\I$ and $\bff{J}_+\restriction_{\partial\Omega_4}-\I$ are strictly upper triangular.
\end{enumerate}
\begin{enumerate}
\item[(iii)]  The matrix $\bff{J}(\lam)$ satisfies the first-order product condition\footnote{The first-order product condition at each intersection point is needed to insure that the RHP has a continuous solution.  See \cite[Definition 2.55]{TO16} for the statement of this condition and 
\cite[Section 2.7]{TO16} for further discussion.}
at the 
\bluetext{ intersection} points $\pm S_\infty$ with the real $\lam$-axis.
\end{enumerate}
\end{theorem}

\end{enumerate}
\end{RHP}
As before, let $\bff{J}_x(\lam) = e^{-ix\lam \ad \sigma_3} \bff{J}$. Solving the Beals-Coifman integral equation 
for 
$$\bff{\nu}=  \bfN_+( \bff{J}_x^+)^{-1}=\bfN_-( \bff{J}_x^-)^{-1},$$
one recovers $q \in H^{2,2}(\R)$ from the relation
\begin{align*}
 q(x) &= 2i\lim_{z \rightarrow\infty} z\, N_{12}(z;x)\\
        &=\left( -\dfrac{1}{\pi}\int_{\Gamma} \bff{\nu}(\lambda;x) 
       e^{-i\lambda x \ad\sigma_3} \left(\bff{J}_+(\lambda) - \bff{J}_- (\lambda) \right)d\lambda \right)_{12}.
\end{align*}

\section{Soliton Resolution}
\label{Resolution}

Theorem~\ref{thm:soliton-resolution} describes the resolution at large times of the solution of \eqref{DNLS} for generic initial data into a sum of a solitonic and radiating components. The rigorous analysis of the large-time behavior of integrable systems with finitely many solitons using RH problems goes back to the pioneering works \cite{DIZ93,DZ93} in which the Deift-Zhou nonlinear steepest descent method was developed. Inspired by previous studies of soliton resolution for other integrable systems \cite{BJM18,CJ16} our analysis uses the $\overline{\partial}$-generalization of the steepest descent method introduced by \cite{DM08, MM06} (see also the recent survey \cite{DMM19}).  As is always the case, we carry out our analysis for \eqref{DNLS2} and, as a last step, use the gauge transformation \eqref{gauge-transform} to obtain results for \eqref{DNLS}.

Our starting point is
\begin{problem}
\label{RHP longtime}
Fix $(x,t) \in \bbR^2$, let $( \rho,  \{ C_i, \lambda_i  \}_{i=1}^{N} )\subset \mathcal{G} \times (\bbC_\times\times \bbC^+)^N$, and let $\Gamma = \mathbb{R} \bigcup_{j=1}^N \left( \Gamma_j \cup \Gamma_j^*\right)$ be oriented as in Figure~\ref{figure-2}.   
Find a vector-valued  function\footnote{In the rest of this section we will often omit the parametric dependence on $x$ and $t$ and write simply $\bfN(z)$.} $\bfN(\dotarg; x,t)$ with the following properties:
\begin{itemize}
\item[(i)](Analyticity) $\bfN(z;x,t )$ is a  row vector-valued  analytic function of $z$ for $z\in \mathbb{C}\setminus \Gamma $.
\item[(ii)] (Normalization) ${\bfN}(z; x,t)=(1 ,0)+\mathcal{O}(z^{-1})$ as $z \rightarrow\infty$.
\item[(iii)] (Jump condition) For each $\lambda\in \Gamma $, $\bfN$ has continuous 
boundary values
$\bfN_{\pm}(\lambda; x,t )$. 
Moreover, the jump relation 
$$\bfN_+(\lambda; x, t)=\bfN_-(\lambda; x,t ) \bff{J}_{x,t}(\lambda)$$ 
holds, where for $\lambda\in\mathbb{R}$
\begin{equation}
\label{N_jump}
\bff{J}_{x,t}(\lambda)=
		e^{-it \theta(\lambda,x/t)\, \ad \sigma_3} 
		\Twomat	{1+\lambda|\rho(\lambda)|^2}{\rho(\lambda)}
						{\lambda\overline{\rho(\lambda)}}{1}
\end{equation}
and
$$
	\theta(z, \xi) = 2 z^2 + \xi z
$$ 
\item[(iv)] (Soliton component)
and for $\lam \in \Gamma_i\cup\Gamma_i^* $ 
$$
\bff{J}_{x,t}(\lambda) = \begin{dcases}
	\twomat{1}{0}{\dfrac{C_i \,\lambda_i e^{2i \theta(\lambda,x/t)}}{\lambda-\lambda_i}}{1}	
	&	\lambda \in \Gamma_i, \\
	\twomat{1}{\dfrac{\overline{C_i} \, e^{-2i \theta( \lambda,x/t)} }{\lambda-\overline{\lambda_i}}}{0}{1}
	&	\lambda \in \Gamma_i^*
\end{dcases}
$$
\end{itemize}
\end{problem}

It  is crucial to the proof of Theorem~\ref{thm:soliton-resolution} that only generic initial data $q_0 \in U$, as defined by Proposition~\ref{thm:R}, are considered, so that the inverse problem under consideration is the finite soliton problem RH problem~\ref{RHP longtime} without spectral singularities on the real line. The absence of spectral singularities immediately implies that there exist positive constants $c_1,c_2>0$ such that the spectral data $( \rho,  \{ C_i, \lambda_i  \}_{i=1}^{N} )$ satisfies
\begin{equation}\label{sp.bounds}
\begin{aligned}
	\frac{1}{c_1} \leq  1 + \lambda | \rho(\lambda)|^2 \leq c_1, \quad \lambda \in \R \\
	 d_\Lambda :=\quad  \ \inf_{\mathclap{\lambda,\mu \in \Lambda_N, \lambda \neq \mu}}  \quad \ \  |\lambda- \mu| > c_2,
	\qquad \Lambda_N = \{ \lambda_k,\lambda_k^*\}_{k=1}^N
\end{aligned}
\end{equation}
Though the inverse spectral problems with arbitrary spectral singularities was introduced in Section~\ref{Arbitrary} for DNLS (see also \cite{BM19,Zhou89-2,Zhou95} for other systems), it is still an open question how to extract asymptotic information from them in the large-time limit. 

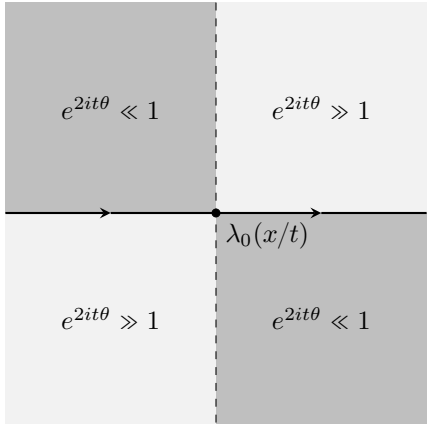
\begin{figure}[htb]
\centering
\begin{tikzpicture}[scale={.7}]
\path[fill=gray!20,opacity=0.5]	(0,0) rectangle(4,4);
\path[fill=gray!20,opacity=0.5]    (-4,-4) rectangle(0,0);
\path[fill=gray,opacity=0.5]		(-4,0) rectangle (0,4);
\path[fill=gray,opacity=0.5]		(0,-4) rectangle (4,0);
\draw[fill] (0,0) circle[radius=0.075];
\draw[thick,->,>=stealth] 	(0,0) -- (2,0);
\draw	[thick]    	(2,0) -- (4,0);
\draw[thick,->,>=stealth]	(-4,0) -- (-2,0);
\draw[thick]		(-2,0) -- (0,0);
\draw[thin,dashed]	(0,4) -- (0,-4);
\node at (2,2) 	{$e^{2it\theta} \gg 1 $};
\node at (-2,-2)  {$e^{2it\theta} \gg 1 $};
\node at (-2,2)	{$e^{2it\theta} \ll 1 $};
\node at (2,-2)	{$e^{2it\theta} \ll 1 $};
\node[below right] at (0,0)		{$\lambda_0(x/t)$};
\end{tikzpicture}
\caption{The regions of growth and decay of the exponential factor $e^{2it \theta(z,x/t)}$ in the complex $z$-plane as $t \to + \infty$.
\label{fig:theta regions}
} 
\end{figure}

The core principal of the Deift-Zhou steepest descent analysis is to use matrix factorizations to introduce transformations, based on factorizations of the jump matrix \eqref{N_jump}, which separate the oscillatory factors $e^{\pm 2it \theta(z,x/t)}$ on the real axis, by moving each onto new contours on which they decay as $t \to +\infty$, see Figure~\ref{fig:theta regions}.  
The useful factorizations for our needs are 
\begin{gather}
	\label{UL}
	\bff{J}_{x,t}(\lambda) = 
	\twomat{1}{\rho(\lambda) e^{-2it \theta(\lambda,x/t)} }{0}{1}
	\twomat{1}{0}{\lambda \overline{\rho(\lambda)} e^{2it \theta(\lambda,x/t)} }{1}, 
\end{gather}
for $\lambda > \lambda_0$ and
\begin{multline}
	\label{LDU}
	\bff{J}_{x,t}(\lambda) = \\[5pt]
	\twomat{1}{0}
	  {\dfrac{\lambda \overline{\rho(\lambda)}}{1+ \lambda |\rho(\lambda)|^2} e^{2it \theta(\lambda,x/t)}} {1} 
	\left( 1+ \lambda |\rho(\lambda)|^2 \right)^{\sigma_3}
	\twomat{1}{\dfrac{\rho(\lambda)}{1+ \lambda |\rho(\lambda)|^2}   e^{-2it \theta(\lambda,x/t)} }
	  {0}{1},
\end{multline}
for $\lam < \lam_0$, 
where 
\begin{equation}\label{crit_pt}
	\lambda_0 = \lambda_0(\xi) := -\frac{\xi}{4}
\end{equation} 
is the unique critical point of the phase function $\theta(z,\xi)$. These factorizations algebraically separate the dependence on the exponential factors. If the reflection coefficient $\rho(\lambda)$ is analytic in a strip containing the real axis, then the right-most (resp. left-most) factor in each factorization has an analytic continuation into the upper (resp. lower) half-plane from the indicated real half-line such that these continuations are exponentially near identity as $t \to \infty$ away from the simple critical point $\lambda_0(x/t)$. For general data, however, there is no reason to suspect that $\rho$ should have an analytic continuation from the real axis, and delicate approximation procedures were introduced \cite{DZ03} in order to consider this data. The innovation of \cite{DM08} was to admit \emph{non-analytic} continuation of $\rho(\lambda)$ by considering the $\overline{\partial}$-generalization of Riemann-Hilbert problems. Each of the above references considers soliton-free inverse problems. In \cite{BJM18} it was demonstrated how to adapt the $\overline{\partial}$-method to studying the large-time asymptotic behavior of a  problem with solitons in such a way that the soliton and dispersive components of the problem could be treated essentially-independently of one another. 

Before describing the steepest descent method in detail below, we summarize our procedure, which can be divided into four steps.
\begin{enumerate}
	\item Make a conjugating transformation $\mathbf{N}^{(1)}(z) = \mathbf{N}(z) \delta(z; \lambda_0)^{-\sigma_3}$, where the function $\delta(z,\lambda_0)$ is a partial transmission coefficient. Algebraically, this transformation has the effect of removing the central diagonal factor from the factorization \eqref{LDU}. 
	
	\smallskip
	
	\item Introduce continuous non-analytic interpolants of the outer factors in \eqref{UL} and \eqref{LDU} from their half-lines to the appropriate sectors of the complex plane.
	
	\smallskip
	
	\item Construct  a global approximate \emph{matrix} solution, which captures the leading order behavior of the solution to the RH problem. The approximation consists of an exterior model away from the critical point which captures the soliton component of the solution and a local model near $\lambda_0$ which captures the dispersive component. 
	
	\smallskip
	
	\item Make rigorous estimates of the residual error using the representation of the solution in terms of the Cauchy singular integrals over $\C$. Having solid Cauchy integrals, bounded $L^\infty (\mathbb{C})$ operators, instead of Cauchy integrals over contours, which is an  $L^2$ theory, is the essential reason that the $\overline{\partial}$ theory is so much simpler than the standard steepest descent method in \cite{DZ03}. 
\end{enumerate}	

\subsection{Step 1: Conjugation}	
Recognizing that we will later use the factorizations \eqref{UL}-\eqref{LDU} to deform the outermost factors in each of these factorizations into regions of decay, we first introduce a diagonal conjugation which has the effect of removing the diagonal central factor in \eqref{LDU} from the half-line $\lambda < \lambda_0$. Let
\begin{equation}\label{delta}
	\delta(z) = \delta(z;\lambda_0) := \exp \left( i \int_{-\infty}^{\lambda_0} \frac{ \kappa(\lambda)}{\lambda-z} d\lambda \right),
	\quad
	\kappa(\lambda) = -\frac{1}{2\pi} \log(1 + \lambda |\rho(\lambda)|^2).
\end{equation}
This function has several useful properties which we summarize in the following proposition, the proof of which is given in \cite[Appendix A]{LPS18a}. 
\begin{lemma}\label{lem:delta}
The function $\delta(z)$ defined by \eqref{delta} has the following properties:
\begin{enumerate}
	\item[(i)] $\delta(z)$ is meromorphic in $\C \setminus (-\infty, \lambda_0 ]$. 

	\item[(ii)] For $z \in \C \setminus (-\infty, \lambda_0 ]$, 
	$\delta(z) \overline{\delta(\overline{z})} = 1$.
	Moreover,
	$\displaystyle
		e^{-\| \kappa \|_\infty/2} \leq 
		\left| \delta(z) \right|
		\leq e^{\| \kappa \|_\infty/2} ~.
	$
	\item[(iii)] For $\lambda \in (-\infty, \lambda_0 ]$, $\delta$'s boundary values  $\delta_\pm$, as $\lambda$ approaches the real axis from above and below, satisfy 
	\begin{equation}\label{delta.jump}
		\delta_+(\lambda) =  \delta_-(\lambda) ( 1 + \lambda |\rho(\lambda)|^2).
	\end{equation}
		\item[(iv)] As $ |z | \to \infty$ with $|\arg(z)| \neq \pi$, 
	\begin{equation}\label{delta.expand}
		{  \delta(\lambda)} = 1 + 
		\frac{\delta_1}{z} 
		+ \mathcal{O} \left( z^{-2} \right),
		\quad \delta_1 = -i \int_{-\infty}^{\lambda_0} \kappa(s) ds .
	\end{equation}
	\item[(v)] As $z \to \lambda_0$ along any ray $\lambda_0 + e^{i \phi} \R_+$ with 
	$| \arg (z-\lambda_0) | < \pi$
	\begin{gather} 
	\nonumber
		\left| \delta(z; \lambda_0) - \delta_0(\lambda_0) (z-\lambda_0)^{i \kappa(\lambda_0)} \right| 
		\lesssim_{\rho,\phi} - |z - \lambda_0 | \log |z - \lambda_0|.
\intertext{The implied constant depends on the $H^{2,2}(\R)$-norm of $\rho$ and is independent of $\lambda_0$. Here $\delta_0(\lambda_0)$ is the complex unit}
		\label{delta0 arg}
		\delta_0(\lambda_0) = e^{i \beta(\lambda_0)},
		\quad 
		\beta(\lambda_0) :=
		 \frac{1}{2\pi} \int_{-\infty}^{\lambda_0} \mathrm{d}_\lambda \log(1+\lambda |\rho(\lambda)|^2 \log(\lambda_0 - \lambda)  d\lambda,
	\end{gather}
In all of the above formulas, we choose the principal branch of power and logarithm functions.
\end{enumerate}
\end{lemma}

Using $\delta(z)$ we define the new unknown 
\begin{equation}\label{N1}
	\mathbf{N}^{(1)}(z) = \begin{dcases}
		\mathbf{N}(z) 
		\Twomat{1}{0}{-C_k \lambda_k \frac{ 1 - (\delta(z)/\delta(z_k))^2}{z-z_k} e^{2it \theta(z,x/t)} }{1}
		\delta(z)^{-\sigma_3} & z \in \operatorname{int}(\Gamma_j) \\[5pt]
		\mathbf{N}(z) 
		\Twomat{1}{-\overline{C_k} \frac{ 1 - (\delta(z)/\delta(z_k))^{-2}}{z-\overline{z_k}} e^{-2it \theta(z,x/t)} }{0}{1}
		\delta(z)^{-\sigma_3} & z \in \operatorname{int}(\Gamma_j^*) \\[5pt]
		\mathbf{N}(z) \delta(z)^{-\sigma_3} & \text{elsewhere}.
	\end{dcases}
\end{equation}
A direct computation\footnote{We observe that as $\delta(z)$ is analytic off $\mathbb{R}$, the off diagonal entries in the middle factor of \eqref{N1} are analytic} shows that the new unknown $\mathbf{N}^{(1)}$ satisfies the following problem
%
%
\begin{RHP}\label{rhp:N1}
Find a vector-valued  function $\bfN^{(1)}(\dotarg)$ with the following properties:
\begin{itemize}
\item[(i)](Analyticity) $\bfN^{(1)}(z)$ is a  row vector-valued  analytic function of $z$ for $z\in \mathbb{C}\setminus \Gamma $.
\item[(ii)] (Normalization) ${\bfN}^{(1)}(z)=(1 ,0)+\mathcal{O}(z^{-1})$ as $z \rightarrow\infty$.
\item[(iii)] (Jump condition) For each $\lambda\in \Gamma $, $\bfN^{(1)}$ has continuous 
boundary values
$\bfN^{(1)}_{\pm}(\lambda)$. 
Moreover, the jump relation 
$$\bfN^{(1)}_+(\lambda)=\bfN^{(1)}_-(\lambda) \bff{J}^{(1)}_{x,t}(\lambda)$$ 
holds, where for $\lambda\in\mathbb{R}$
\begin{gather}
\label{N1_jump}
\bff{J}^{(1)}_{x,t}(\lambda)= 
   \begin{dcases}
       \Twomat{1}{ \rho(\lambda) e^{-2it \theta(\lambda,x/t)} \delta(\lambda)^2 }{0}{1}
       \Twomat{1}{0}{ \lambda \overline{\rho(\lambda)} e^{2it \theta(\lambda,x/t)} \delta(\lambda)^{-2} }{1}
	& \lambda > \lambda_0 \\
       \Twomat{1}{0}{ \frac{\lambda \overline{\rho(\lambda)} \delta_-(\lambda)^{-2} }{1+ \lambda |\rho(\lambda)|^2} e^{2it \theta(\lambda,x/t)}  }{1}
       \Twomat{1}{ \frac{\rho(\lambda)  \delta_+(\lambda)^{2} }{1+ \lambda |\rho(\lambda)|^2} e^{-2it \theta(\lambda,x/t)}}{0}{1}
	& \lambda < \lambda_0
   \end{dcases}
\end{gather}

\item[(iv)] (Soliton component)
For $\lam \in \Gamma_k\cup\Gamma_k^* $ 
$$
\bff{J}_{x,t}(\lambda) = \begin{dcases}
	\twomat{1}{0}{\dfrac{C_k \delta(\lambda_k)^{-2} \lambda_k }{\lambda-\lambda_k}  e^{2i \theta(\lambda,x/t)} }{1}	
	&	\lambda \in \Gamma_k, \\
	\twomat{1}{\dfrac{\overline{C_k} \delta(\lambda_k)^{2} }{\lambda-\overline{\lambda_k}} e^{-2i \theta( \lambda,x/t)}}{0}{1}
	&	\lambda \in \Gamma_k^*
\end{dcases}
$$
\end{itemize}
\end{RHP}

\subsection{Step 2: Non-analytic extensions}
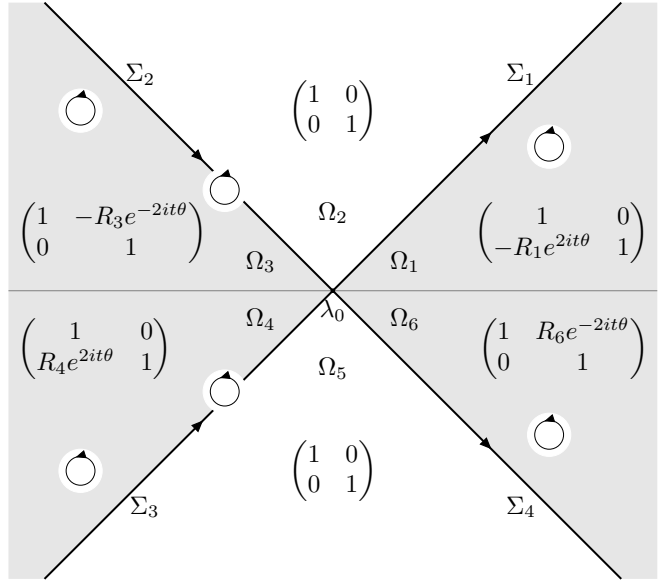
\begin{figure}[h]
\centering
\resizebox{0.7\textwidth}{!}{
\begin{tikzpicture}
\path [fill=gray!20] (0,0) -- (-4,4) -- (-4.5,4) -- (-4.5,-4) -- (-4,-4) -- (0,0);
\path [fill=gray!20] (0,0) -- (4,4) -- (4.5,4) -- (4.5,-4) -- (4,-4) -- (0,0);
%
\draw [help lines] (-4.5,0) -- (4.5,0);
\draw [thick][->-] (-4,4) -- (0,0);
\draw [thick][->-] (-4,-4) -- (0,0);
\draw [thick][->-] (0,0) -- (4,4);
\draw [thick][->-] (0,0) -- (4,-4);
%
\foreach \pos in { 
		(-1.5,1.4), (-1.5,-1.4), (-3.5,2.5), 
		(-3.5,-2.5), (3,2), (3,-2)}
{	
\draw[color=white, fill=white] \pos circle [radius=.3];
\draw[color=black][->-] \pos++(0,-.2) arc (-90:270:0.2);
}
%
\node[left] at (3,3) {$\Sigma_{1}\,$};
\node[right] at (-3,3) {$\Sigma_{2}$};
\node[right] at (-3,-3) {$\,\Sigma_{3}$};
\node[left] at (3,-3) {$\Sigma_{4}\,$};
\draw[fill] (0,0) circle [radius=0.025];
\node[below] at (0,0) {$\lambda_0$};
%
\node at (1,.4) {$\Omega_{1}$};
\node at (0,1.08) {$\Omega_{2}$};
\node at (-1,.4) {$\Omega_{3}$};
\node at (-1,-.4) {$\Omega_{4}$};
\node at (0,-1.08) {$\Omega_{5}$};
\node at (1,-.4) {$\Omega_{6}$};
%
\node[left] at (4.5,0.8) {$\twomat{1}{0}{-R_1 e^{2it \theta}}{1}  $ };
\node[right] at (-4.5,0.8) {$\twomat{1}{-R_3 e^{-2it \theta}}{0}{1} $ };
\node[right] at (-4.5,-0.8) {$\twomat{1}{0}{R_4 e^{2it \theta}}{1} $ };
\node[left] at (4.5,-0.8) {$\twomat{1}{ R_6 e^{-2it \theta }}{0}{1} $};
\node at (0,2.5) {$\twomat{1}{0}{0}{1}$} ;
\node at (0,-2.5){$\twomat{1}{0}{0}{1}$} ;
\end{tikzpicture}
}
\caption{
The contour $\Gamma^{(2)}$ and sectors $\Omega_k,\ k=1,\dots, 6$, in the $z = u + i v$-plane. Shaded regions indicate the support of the $\dbar$ derivative of the non-analytic extensions defining the transformation from $\mathbf{N}^{(1)} \mapsto \mathbf{N}^{(2)}$.
\label{fig:extensions}
}
\end{figure}

We're now ready for the key step in the steepest descent analysis, 
moving the oscillatory factors $e^{\pm it \theta(z,\lambda_0) }$ from the real axis onto appropriate contours in the complex plane where they decay as $t \to +\infty$. 
To do this, we introduce non-analytic extensions of the coefficients of $e^{\pm i \theta(z,\lambda_0) }$ in each of the factors in \eqref{N1_jump}. As depicted in Figure~\ref{fig:extensions}, we define the contours, 
\begin{gather}\label{Lambda.2}
	\Gamma^{(2)} = 
	\Sigma \bigcup_{k=1}^N ( \Gamma_k \cup \Gamma_k^*) \\
	\nonumber
	\Sigma = \bigcup_{j=1}^4 \Sigma_k, \quad 
	\Sigma_k = \lambda_0 + e^{ \frac{i\pi}{4} (2k-1)} \mathbb{R}_+, \ k=1,2,3,4.  
\end{gather}
where each $\Sigma_k$ is oriented with increasing real part. We denote the six connected components of $\mathbb{C} \setminus \left( \mathbb{R} \cup  \Sigma \right)$ by $\Omega_k$, $k=1,\dots, 6$, numbered counterclockwise starting from the half-line $z > \lambda_0$.

Writing $z = u + i v$ and letting $\Lambda_N$ denote the discrete spectral data as in \eqref{sp.bounds}, we choose $\chi_{\Lambda}(u,v)$, a $C_0^\infty(\mathbb{C},[0,1])$ cutoff function supported on a neighborhood of each point of the discrete spectrum such that
\[
	\chi_\Lambda(u,v) = \begin{cases}
		1 & \dist(u+iv, \Lambda_N)  < d_\Lambda/3 \\
		0 & \dist(u+iv, \Lambda_N)  > 2d_\Lambda/3 
	\end{cases}
\]
so that the disks of support of $\chi_\Lambda$ intersect neither each other nor the real axis. We also choose the contours $\Gamma_j \cup \Gamma_j^*$ such that they lie entirely in the set $\dist(z, \Lambda_N) < d_\Lambda/3$ where the $(1-\chi_\Lambda)$ is identically zero. We now define the extensions, supported on the shaded domains in Figure~\ref{fig:extensions}, 
as follows.
\begin{subequations}\label{extensions}
\begin{align}
&\text{For } z = u+iv \in \Omega_1 \nonumber \\
&\begin{multlined}[.9\displaywidth]
	R_1(u,v) = \left[ \cos(2\arg(u-\lambda_0 + i v)) z \overline{\rho(u)} \delta(z)^{-2}  
	\right. \\ \left.
	+ \Big( 1 - \cos(2\arg(u-\lambda_0+iv)) \Big) 
	\lambda_0 \overline{\rho(\lambda_0)} \delta_0(\lambda_0)^{-2} (z-\lambda_0)^{-2i \kappa(\lambda_0)} 
	\right] (1 - \chi_{\Lambda}(u,v))
\end{multlined} \\
&\text{For } z = u+iv \in \Omega_3 \nonumber \\
&\begin{multlined}[\displaywidth]
	R_3(u,v) = \left[ \cos(2\arg(u-\lambda_0 + i v))  \frac{\rho(u)}{1+ z |\rho(u)|^2} \delta(z)^{2} 
	\right. \\ \left.
	+ \Big( 1 - \cos(2\arg(u-\lambda_0+iv)) \Big)
	\frac{\rho(\lambda_0)}{1+\lambda_0 |\rho(\lambda_0)|^2} 
	\delta_0(\lambda_0)^{2} (z-\lambda_0)^{2i \kappa(\lambda_0)} 
	 \right] (1 - \chi_{\Lambda}(u,v))
\end{multlined} \\
&\text{For } z = u+iv \in \Omega_4 \nonumber \\
&\begin{multlined}[\displaywidth]
	R_4(u,v) = \left[ \cos(2\arg(u-\lambda_0 + i v)) \frac{z\overline{\rho(u)}}{1+ z |\rho(u)|^2} \delta(z)^{-2}  
	\right. \\ \left.
	+ \Big( 1 - \cos(2\arg(u-\lambda_0+iv)) \Big) \frac{\lambda_0\overline{\rho(\lambda_0)}}{1+\lambda_0 |\rho(\lambda_0)|^2} 
	\delta_0(\lambda_0)^{-2} (z-\lambda_0)^{-2i \kappa(\lambda_0)} 
	 \right] (1 - \chi_{\Lambda}(u,v))
\end{multlined} \\
&\text{For } z = u+iv \in \Omega_6 \nonumber \\
&\begin{multlined}[\displaywidth]
	R_6(u,v) = \left[ \cos(2\arg(u-\lambda_0 + i v))  \rho(u) \delta(z)^{2} 
	\right. \\ \left.
	+ \Big( 1 - \cos(2\arg(u-\lambda_0+iv)) \Big) \rho(\lambda_0) \delta_0(\lambda_0)^{2} (z-\lambda_0)^{2i \kappa(\lambda_0)} 
	\right] (1 - \chi_{\Lambda}(u,v))
\end{multlined}
\end{align}
\end{subequations}

It's straightforward to check that each extension $R_k(u,v)$ is continuous as $\Omega_k \ni (u,v) \to (\lambda_0,0)$ and that: 
\begin{itemize}
	\item $R_1(u,v) = z \overline{\rho(z)} \delta(z)^{-2}$ for $z =u+iv \in (\lambda_0,\infty)$;  
	\smallskip
	\item $R_3(u,v) = \frac{\rho(z)}{1+z |\rho(z)|^2} \delta(z)^{2}$ for $z =u+iv \in (-\infty, \lambda_0)$;  
	\smallskip
	\item $R_4(u,v) = \frac{z \overline{\rho(z)}}{1+z |\rho(z)|^2} \delta(z)^{-2}$ for $z =u+iv \in (-\infty, \lambda_0)$; 
	\smallskip
	\item $R_6(u,v) = \rho(z) \delta(z)^{2}$ for $z =u+iv \in (\lambda_0,\infty)$. 
\end{itemize}
Moreover, along the diagonal boundary of each sector the extensions evaluate to simple constants times a factor of the form $(z - \lambda_0)^{\pm 2 i \kappa(\lambda_0)}$ which, as we will see, will be very convenient for introducing a local model at the critical point $\lambda_0$. 

We now ``open lenses'' off each real half-lines $\lambda \lessgtr \lambda_0$, by using the extensions defined above to introduce a change of variables which trades the oscillator jumps along the real axis for new jumps along the contours $\Sigma_k$ where the exponential factors $e^{\pm 2 i t \theta(z;x/t)}$ are exponentially decaying as $t \to + \infty$. Our change of variables is:  

\begin{equation}\label{N2}
\begin{gathered}
	\mathbf{N}^{(2)}(u,v;x,t) = \mathbf{N}^{(1)}(u+iv; x,t) \mathbf{R}(u,v;x,t) \\[5pt]
	\mathbf{R}(u,v;x,t) = \begin{dcases}
	\twomat{1}{0}{R_1(u,v) e^{2it \theta(z,x/t)} }{1}^{-1} 
	& z = u+iv \in \Omega_1 \\
	\\
	\twomat{1}{R_3(u,v) e^{-2it \theta(z,x/t)} }{0}{1}^{-1} 
	& z = u+iv \in \Omega_3 \\
	\\
	\twomat{1}{0}{R_4(u,v) e^{2it \theta(z,x/t)} }{1}
	& z = u+iv \in \Omega_4 \\
	\\
	\twomat{1}{R_6(u,v) e^{-2it \theta(z,x/t)} }{0}{1}
	& z = u+iv \in \Omega_6 \\
	\\
	\I  & z = u+iv \in \Omega_2 \cup \Omega_5
\end{dcases}
\end{gathered}
\end{equation}
The notation $\mathbf{N}^{(2)}(u,v;x,t)$ reflects the fact that, because of the non-analytic extensions $R_j$, $j=1,3,4,6$, $\mathbf{N}^{(2)}$ is not a piecewise-analytic function. However, as discussed above, the extensions $R_j$, $j=1,3,4,6$, are continuous in the relevant sectors with continuous extension to their boundary, so the function $\mathbf{N}^{(2)}(u,v;x,t)$ is a piecewise-continuous function of $(u,v) \in \mathbb{R}^2$ with jump discontinuities along the sector boundaries and the disks  boundaries $\bigcup_{k=1}^n (\Gamma_k \cup \Gamma_k^*)$. 

The following Lemma will be  needed later when we study the $\dbar$-problem that results from the nonanalytic transformation \eqref{N2}.
\begin{lemma}[\cite{JLPS18b} Lemma~2.4] 
\label{lem:R.dbar}
Suppose that $\rho \in H^{2,2}(\R)$ and that $$c = \inf_{\lambda \in \R} (1 + \lambda |\rho(\lambda)|^2) $$ is strictly positive. Then for $u+iv \in \Omega_k,\ k=1,3,4,6$, we have
the estimates
$$
\left| \dbar \mathbf{R}(u,v) \right| 
	\lesssim 
	(| \dbar \chi_\Lambda(u,v)| + | R_{k}(u,v) |  - \log | u + i v -\lambda_0 |)
	  e^{-8t |uv |}
$$
for $ |u+ iv - \lambda_0| \leq 1$,
$$
\left| \dbar \mathbf{R}(u,v) \right| 
	\lesssim
	 (| \dbar \chi_\Lambda(u,v)| + | R_{k}(u,v) |(1+| u+iv -\lambda_0|^2)^{-1/2} )
	    e^{-8t |uv |}
$$
for $ |u+ iv - \lambda_0| > 1	$,
and 
$$
\dbar \mathbf{R}(u,v) = 0 
$$
if $u+iv \in \Omega_2 \cup \Omega_5$.	
Here, $R_{k}(u,v)$ represents the nonzero off-diagonal entries of $\mathbf{R}$ for $u+iv \in \Omega_k,\ k=1,3,4,6$ defined by \eqref{extensions} and \eqref{N2}.
\end{lemma}

\subsection{Step 3: Constructing a matrix model solution}
Recognizing that the exponential factors $e^{\pm 2i t \theta(z,x/t)}$ in \eqref{N2} decay exponentially to zero as $t \to +\infty$ pointwise within the indicated sectors, it's reasonable to suppose that the contribution of the non-analytic $\overline{\partial}$ derivatives of $\mathbf{N}^{(2)}$ are negligible at large times. 

Simply ignoring the $\overline{\partial}$ derivative of $\mathbf{N}^{(2)}$ we arrive at a Riemann-Hilbert problem for a matrix function $\mathbf{P}(z;x,t)$ whose jumps are exactly equal to those of $\mathbf{N}^{(2)}(\lambda;x,t)$ along $\Gamma^{(2)}$ (cf. \eqref{Lambda.2}).
We introduce, for convenience, the scaled local variable and phase function
\begin{gather}\label{pc.variables}
	\zeta = \zeta(z;x,t) := (8t)^{1/2} (z-\lambda_0) \\
	\omega(x,t) := \arg \rho(\lambda_0) + 2 \beta(\lambda_0) 
	- \kappa(\lambda_0) \log(8t) + 4t \lambda_0^2,
\end{gather}
where $\kappa(\lambda_0)$ and $\beta(\lambda_0)$ are defined by \eqref{delta} and \eqref{delta0 arg} respectively.
Then $\mathbf{P}(z;x,t)$ solves the following Riemann-Hilbert problem.
\begin{problem}\label{model}
Fix $x \in \R$ and $t>0$. Find a $2\times2$ matrix-valued function $\mathbf{P}(\dotarg;x,t)$ with the following properties:
\begin{itemize}
\item[(i)](Analyticity) $\mathbf{P}(z;x,t)$ an analytic function of $z$ for $z\in \mathbb{C}\setminus\Gamma^{(2)} $.
\item[(ii)](Symmetry) The entries of $\mathbf{P}(z;x,t)$ satisfies the relations 
\[
\mathbf{P}_{22}(z;x,t) = \overline{\mathbf{P}_{11}(\overline{z};x,t)}
\qquad
\mathbf{P}_{21}(z;x,t) = -z \overline{\mathbf{P}_{12}(\overline{z};x,t)}
\]
\item[(iii)] (Normalization) $\mathbf{P}(z;x,t)= \begin{psmallmatrix} 1 &0 \\ \alpha & 1 \end{psmallmatrix} +\mathcal{O}(z^{-1})$ as $z \rightarrow\infty$, for a constant $\alpha$ determined by the symmetry condition above.
\item[(iv)] (Jump condition) For each $\lambda \in \Gamma^{(2)}$, the boundary values $\mathbf{P}_+(z;x,t)$ (resp. $\mathbf{P}_-(z;x,t)$), taken along any component of $\Gamma^{(2)}$ from the left (resp. right) according to the orientation of $\Gamma^{(2)}$ shown in Figure~\ref{fig:extensions}, 
are related  by
$\mathbf{P}_+(\lambda;x,t)=\mathbf{P}_-(\lambda;x,t) \mathbf{J}^{(P)}_{x,t}(\lambda)$
where
\begin{gather}
\label{P.jump}
\mathbf{J}^{(P)}_{x,t}(\lambda)= 
   \begin{dcases}
       \I + (1 - \chi_{\Lambda}(u,v)) 
       \twomat{0}{0}{ \lambda_0 |\rho(\lambda_0)| e^{-i \omega(x,t)} \zeta^{-2i \kappa(\lambda_0)} e^{i\zeta^2/2} }{0}
	& \lambda \in \Sigma_1  \\
      \I + (1 - \chi_{\Lambda}(u,v)) 
      \twomat{0}{\frac{|\rho(\lambda_0)|}{1+\lambda_0 |\rho(\lambda_0)|^2} 
	e^{i \omega(x,t)} \zeta^{2i \kappa(\lambda_0)} e^{-i\zeta^2/2} } {0} {0}
       & \lambda \in \Sigma_2 \\
       \I + (1 - \chi_{\Lambda}(u,v)) 
       \twomat{0}{0} {\frac{ \lambda_0 |\rho(\lambda_0)|}{1+\lambda_0 |\rho(\lambda_0)|^2} 
	e^{-i \omega(x,t)} \zeta^{-2i \kappa(\lambda_0)} e^{i\zeta^2/2} } {0}
       & \lambda \in \Sigma_3 \\
       \I + (1 - \chi_{\Lambda}(u,v))
       \twomat{0}{ |\rho(\lambda_0)| e^{i \omega(x,t)} \zeta^{2i \kappa(\lambda_0)} e^{-i\zeta^2/2}  }{0}{0}
	& \lambda \in \Sigma_4
   \end{dcases}
\end{gather}

\item[(iv)] (Soliton component)
and for $\lam \in \Gamma_k\cup\Gamma_k^* $ 
\begin{equation}\label{P.soliton}
\mathbf{J}_{x,t}^{(P)}(\lambda) = \begin{dcases}
	\twomat{1}{0}{\dfrac{\widetilde{C}_k(x/t) \lambda_k }{\lambda-\lambda_k}  e^{2i \theta(\lambda,x/t)} }{1}	
	&	\lambda \in \Gamma_k, \\
	\twomat{1}{\dfrac{\overline{\widetilde{C}_k(x/t)}  }{\lambda-\overline{\lambda_k}}   e^{-2i \theta( \lambda,x/t)}}{0}{1}
	&	\lambda \in \Gamma_k^*
\end{dcases}
\end{equation}
\[
	\widetilde{C}_k(\xi) = C_k \delta(z_k;\lambda_0(\xi))^{-2} 
	= C_k \exp \left( \frac{i}{\pi} \int_{-\infty}^{\lambda_0(\xi)} \log(1 + s | \rho(s)|^2) \frac{ds}{s-\lambda_k} \right).
\]
\end{itemize}
\end{problem}	

This problem is ideally suited to asymptotic analysis as the jumps \eqref{P.jump}-\eqref{P.soliton} have the following properties:
\begin{itemize}
	\item Outside any open neighborhood of $\lambda_0$ the jumps \eqref{P.jump} are uniformly exponentially near the identity matrix as $t \to \infty$, suggesting that a good approximation is to consider only the jump \eqref{P.soliton}.
	\item Inside a small neighborhood of $\lambda_0$, bounded away from the support of \eqref{P.soliton}, the jumps \eqref{P.jump} match exactly a well-known problem in the literature, the so-called parabolic cylinder model, which can be solved exactly. 
\end{itemize}
Define the neighborhood of the critical point
\[
	\mathcal{U}_{\lambda_0} = \{ z \in \mathbb{C}\, : \, |z-\lambda_0| < \delta_\Lambda/3\}
\] 
where the radius is chosen such that  so as to be bounded away from the soliton contours $\bigcup_{k=1}^n (\Gamma_k \cup \Gamma_k^*)$. As we outline below we construct the solution $\mathbf{P}$ of Problem~\ref{model} by patching together two explicit approximates  $\mathbf{P}^{\mathrm{pc}}$ and $\mathbf{P}^{\mathrm{sol}}$ inside and outside $\mathcal{U}_{\lambda_0}$ (respectively) as follows:
\begin{equation}\label{model.form}
	\mathbf{P}(z;x,t) = \begin{cases}
		\mathbf{E}(z;x,t) \mathbf{P}^{\mathrm{sol}}(z;x,t) & \lambda \notin \mathcal{U}_{\lambda_0} \\
		\mathbf{E}(z;x,t) \mathbf{P}^{\mathrm{sol}}(z;x,t) \mathbf{P}^{\mathrm{pc}}(z;x,t) 
		& \lambda \in \mathcal{U}_{\lambda_0}.
	\end{cases}
\end{equation}
This equation serves an implicit definition of the residual error $\mathbf{E}(z;x,t)$. Writing down the resulting Riemann-Hilbert problem for $\mathbf{E}$ one shows that its jump matrix $\mathbf{J}^E$ is uniformly asymptotically near the identity. For such near-identity problems there is a well established existence theory that allows one to compute an asymptotic expansion for $\mathbf{E}$ which complete the construction of $\mathbf{P}$.

\subsubsection*{The soliton model} Let $\mathbf{P}^\mathrm{sol}$ solve the following RHP.
%
%
\begin{RHP}\label{soliton model}
Take $x \in \R$ and $t>0$ with $\xi = x/t$ fixed. Find a $2\times2$ matrix-valued function $\mathbf{P}^{\mathrm{sol}}(\dotarg;x,t)$ with the following properties:
\begin{itemize}
\item[(i)](Analyticity) $\mathbf{P}^{\mathrm{sol}}(z;x,t)$ an analytic function of $z$ for $z\in \mathbb{C}\setminus \bigcup_{k=1}^N (\Gamma_k \cup \Gamma_k^*) $.
\item[(ii)](Symmetry) $\mathbf{P}^{\mathrm{sol}}(z;x,t)$ satisfies the relations 
${P}^{\,\mathrm{sol}}_{22}(z;x,t) = \overline{{P}^{\,\mathrm{sol}}_{11}(\overline{z};x,t)}$ and
${P}^{\,\mathrm{sol}}_{21}(z;x,t) = -z \overline{{P}^{\,\mathrm{sol}}_{12}(\overline{z};x,t)}$.
\item[(iii)] (Normalization) $\mathbf{P}^{\mathrm{sol}}(z;x,t)= \begin{psmallmatrix} 1 &0 \\ \alpha & 1 \end{psmallmatrix} +\mathcal{O}(z^{-1})$ as $z \rightarrow\infty$, for a constant $\alpha$ determined by the symmetry condition above.
\item[(iv)] (Jump condition) For $\lambda \in \bigcup_{k=1}^N (\Gamma_k \cup \Gamma_k^*)$ the boundary values of $\mathbf{P}^{\mathrm{sol}}(z;x,t)$ satisfy the jump relation
$\mathbf{P}^{\mathrm{sol}}_+(\lambda;x,t)=\mathbf{P}^{\mathrm{sol}}_-(\lambda;x,t) \mathbf{J}^{(P)}_{x,t}(\lambda)$
where $\mathbf{J}^{\mathrm{sol}}_{x,t}(\lambda)$ is equal to $\mathbf{J}^{(P)}$ (cf. \eqref{P.soliton}) restricted to $\bigcup_{k=1}^N (\Gamma_k \cup \Gamma_k^*)$.
\end{itemize}
\end{RHP}	

The following proposition is proved in \cite{JLPS18b}:
\begin{proposition}\label{prop:Nsol}
Given scattering data $(\rho, \{C_k, \lambda_k\}_{k=1}^N) \subset H^{2,2}(\mathbb{R} \times (\mathbb{C}_\times \times \mathbb{C}^+)^N$, take $x \in \mathbb{R}$ and $t>0$ with $\xi = x/t$ fixed. Then there exists a unique solution $\mathbf{P}^{\textrm{sol}}(z;x,t)$ of RH Problem~\ref{soliton model} which is exactly the solution of RH Problem~\ref{RHP-2} given the reflectionless scattering data $\mathcal{D}_\xi= \{\widetilde{C}_k, \lambda_k\}_{k=1}^N $ generated by an $N$-soliton solution $q_\mathrm{sol}(x,t; \mathcal{D}_\xi)$ of \eqref{DNLS2} where $\widetilde{C}_k$ is related to ${C}_k$ by \eqref{P.soliton}. Moreover, $\det \mathbf{P}^\mathrm{sol}(z;x,t) = 1$ and  
\[
	\| \mathbf{P}^{\mathrm{sol}}(\dotarg;x,t)  \|_{L^\infty (\mathbb{C})} = \mathcal{O}(1) 
\]
where the implied constant is independent of $(x,t)$ and depends on $\rho$ through its $H^{2,2}(\mathbb{R})$ norm.
\end{proposition}

\subsubsection*{The parabolic cylinder model}
Inside $\mathcal{U}_{\lambda_0}$ we ignore the soliton component jumps \eqref{P.soliton} to arrive at the following local model:

\begin{problem}\label{PCmodel}
Fix $x \in \R$ and $t>0$. Find a $2\times2$ matrix-valued function $\mathbf{P}^{\mathrm{pc}}(\dotarg;x,t)$ with the following properties:
\begin{itemize}
\item[(i)](Analyticity) $\mathbf{P}^{\mathrm{pc}}(z;x,t)$ an analytic function of $z$ for $z\in \mathbb{C}\setminus \ \Sigma $.
\item[(ii)] (Normalization) $\mathbf{P}^{\mathrm{pc}}(z;x,t)= \I +\mathcal{O}(z^{-1})$ as $z \rightarrow\infty$.
\item[(iv)] (Jump condition) For $\lambda \in \Sigma$ the boundary values of $\mathbf{P}^{\mathrm{pc}}$ are related  by
$\mathbf{P}^{\mathrm{pc}}_+(\lambda;x,t)=\mathbf{P}^{\mathrm{pc}}_-(\lambda;x,t) \mathbf{J}^{\mathrm{pc}}_{x,t}(\lambda)$
where
\begin{gather}
\label{pc.jump}
\mathbf{J}^{\mathrm{pc}}_{x,t}(\lambda)= 
   \begin{dcases}
       \twomat{1}{0}{ \lambda_0 |\rho(\lambda_0)| e^{-i \omega(x,t)} \zeta^{-2i \kappa(\lambda_0)} e^{i\zeta^2/2} }{1}
	& \lambda \in \Sigma_1  \\
	\\
      \twomat{1}{\frac{|\rho(\lambda_0)|}{1+\lambda_0 |\rho(\lambda_0)|^2} 
	e^{i \omega(x,t)} \zeta^{2i \kappa(\lambda_0)} e^{-i\zeta^2/2} } {0} {1}
       & \lambda \in \Sigma_2 \\
      \\
       \twomat{1}{0} {\frac{ \lambda_0 |\rho(\lambda_0)|}{1+\lambda_0 |\rho(\lambda_0)|^2} 
	e^{-i \omega(x,t)} \zeta^{-2i \kappa(\lambda_0)} e^{i\zeta^2/2} } {1}
       & \lambda \in \Sigma_3 \\
      \\
       \twomat{1}{ |\rho(\lambda_0)| e^{i \omega(x,t)} \zeta^{2i \kappa(\lambda_0)} e^{-i\zeta^2/2}  }{0}{1}
	& \lambda \in \Sigma_4
   \end{dcases}
\end{gather}
\end{itemize}
\end{problem}

The jump matrices of Problem \eqref{PCmodel} are equal to those of Problem \eqref{model} inside $\mathcal{U}_{\lambda_0}$ and can be solved exactly in terms of special functions, owing to the special choice of $\overline{\partial}$-extensions of the matrix factors in \eqref{extensions}. This model is well-known in the integrable systems literature. The solution of this model is expressed in terms of parabolic cylinder functions, $D_a(z)$ (see \cite[Chapter 12]{DLMF} for its properties); its construction goes back to \cite{Its, DZ93}, see also \cite{LPS18a} for its derivation in the context of DNLS. Here we provide only the necessary formulas 

\begin{proposition}
Let $c_1,c_2$ be strictly positive constants such that $\| \rho \|_{H^{2,2}(\R)} \leq c_1$ and $\inf_{\lambda \in \R} (1+z|\rho(z)|^2) \geq c_2$. Fix $\xi$, let $\kappa = \kappa(\xi)$ be given by \eqref{delta}, and define $\zeta = \zeta(z;x,t)$ and $\omega(x,t)$ as in \eqref{pc.variables}. Then the solution of Problem~\ref{PCmodel} is given by
\[
	\mathbf{P}^{\mathrm{pc}}(z;x,t) = 
	e^{i\omega(x,t) \sigma_3/2} 
	\boldsymbol{\Phi}(\zeta(z)) \mathbf{L}(\zeta(z)) \zeta(z)^{-i \kappa \sigma_3} e^{i\zeta^2 \sigma_3/4} 
	e^{-i\omega(x,t) \sigma_3/2}
\]
where
\begin{gather}
\nonumber
	\mathbf{L}(\zeta) = 
	\left\{ 
		\begin{array}{c@{\quad}l@{\hspace{2em} }c@{\quad}l}
			\twomat{1}{0}{\lambda_0 |\rho(\lambda_0)| }{1} 
			& \arg \zeta \in \left( 0, \frac{\pi}{4} \right), &
			\twomat{1}{\frac{ |\rho(\lambda_0)|} { 1+\lambda_0 | \rho(\lambda_0) |^2} }{0}{1} 
			& \arg \zeta \in \left( \frac{3\pi}{4},\pi \right), 
		\medskip \\
			\twomat{1}{- | \rho(\lambda_0)|}{0}{1} 
			& \arg \zeta \in \left( -\frac{\pi}{4}, 0 \right),  &
			\twomat{1}{0}{\frac{-\lambda_0|\rho(\lambda_0)|}{1+\lambda_0 |\rho(\lambda_0)|^2}}{1} 
			& \arg \zeta \in \left( -\pi, -\frac{3\pi}{4} \right), \\ \\
		\multicolumn{4}{c}{
		\I \quad |\arg \zeta| \in \left( \frac{\pi}{4},\frac{3\pi}{4} \right),
		} 
		\end{array}
	\right.
\shortintertext{and defining}
\label{betas}
	\beta_{12} 
	  = \frac{ \sqrt{2\pi} e^{-\pi \kappa/2} e^{i \pi/4} }{\lambda_0|\rho(\lambda_0)| \Gamma(-i\kappa)}
	\qquad
	\beta_{21}
	  = \frac{ -\sqrt{2\pi}e^{-\pi \kappa/2} e^{-i \pi/4} }{|\rho(\lambda_0)| \Gamma(i\kappa)},
\shortintertext{$\boldsymbol{\Phi}$ is given by}
\nonumber
	\Phi (\zeta)= 
	\begin{pmatrix}
 		{e^{-\frac{3\pi}{4}\kappa} D_{i\kappa}(\zeta e^{-3i\pi/4})} &
		-i \beta_{12}
			e^{\frac{\pi}{4}(\kappa-i)} 
		  D_{-i\kappa-1}(\zeta e^{-\pi i/4}) \\[5pt]
		i \beta_{21}
			e^{-\frac{3\pi}{4}(\kappa+i)}
		  D_{i\kappa -1}(\zeta e^{-3i\pi/4}) &
		e^{\pi\kappa/4}D_{-i\kappa}(\zeta e^{-i\pi/4})
	\end{pmatrix}
\shortintertext{for $\Im(\zeta)>0$, and for $\Im(\zeta)<0$}
\nonumber
 	\Phi (\zeta)= 
	\begin{pmatrix}
		e^{\pi\kappa/4}D_{i\kappa}(\zeta e^{\pi i/4}) &
		-i \beta_{12}
			e^{-\frac{3\pi}{4}(\kappa-i)}
				D_{-i\kappa-1}(\zeta e^{3i\pi/4}) \\[5pt]
		i \beta_{21}
			e^{\frac{\pi}{4}(\kappa+i)}
				D_{i\kappa-1}(\zeta e^{\pi i/4}) &
		e^{-3\pi \kappa/4}D_{-i\kappa}(\zeta e^{3i\pi/4})
	\end{pmatrix}.
\end{gather}
Moreover, we have $\det \mathbf{P}^{\mathrm{pc}} (z;x,t) = 1$ and $\| \mathbf{P}^{\mathrm{pc}} (\dotarg;x,t) \|_\infty = \mathcal{O}(1)$ where the implied constant are uniform in $x$ and $t>1$ and depend only on $c_1$ and $c_2$. 
\end{proposition}

The essential fact need about $\mathbf{N}^\mathrm{pc}$, which follows from the rescaled local variable $\zeta(z)$ and the large $\zeta$ behavior of $\boldsymbol{\Phi}(\zeta)$, is that
\[
	\mathbf{P}^\mathrm{pc}(z;x,t) = \I + \frac{(8t)^{-1/2}}{z-\lambda_0} 
	\twomat{0}{-i \beta_{12}e^{-i \omega(x,t)} }{i \beta_{21} e^{i \omega(x,t)} }{0} + \bigO{t^{-1}}, \quad z \in \partial \mathcal{U}_{\lambda_0}.
\]

\subsubsection*{Computing the residual $\mathbf{E}$}
The models $\mathbf{P}^\mathrm{sol}$ and $\mathbf{P}^\mathrm{pc}$ are bounded unimodular functions which exactly match the jumps of $\mathbf{P}$ on the discrete spectral disks $\bigcup_{k=1}^n (\Gamma_k \cup \Gamma_k^*)$ and on the rays $\Sigma_j, j=1, \dots,4$ inside $\mathcal{U}_{\lambda_0}$ respectively. Using \eqref{model.form} as a definition for the residual $\mathbf{E}$ one can show that it solves a Riemann-Hilbert problem with a jump matrix $\mathbf{J}^{E}$ supported on $\Gamma^{E} = \partial \mathcal{U}_{\lambda_0} \cup (\Sigma \setminus \mathcal{U}_{\lambda_0})$ which satisfies
\[
	\| \mathbf{J}^{E} - \I \|_{L^\infty(\Gamma^{E}) \cap L^{2,k}(\Gamma^{E}) }
	= \bigO{ t^{-1/2}}.
\]
Using this estimate and the well-known existence theorem for Riemann-Hilbert problems with near identity jump matrices \cite{DIZ93,TO16,Zhou89-1} one can express the solution $\mathbf{E}$ in terms of the Cauchy integral operator
\[
	C_{E} f = C^-( f (\mathbf{J}^{E} - \I) ), 
	\qquad 
	C^{-}f (\lambda) = \lim_{\lambda \to \Gamma^{E}_- } \frac{f(z) } {z-\lambda} dz. 
\]
Using the fact that the Cauchy projection operator $C^-$ is bounded one uses the bounds on the jump matrix to show that the function $\mathbf{E}$ exists and compute the asymptotic expansion of the solution for $t \gg 1$. The details of this calculation for DNLS are found in \cite[Section 3.3]{JLPS18b}. The result of that computation gives:
\begin{equation}\label{model.error}
  \begin{gathered}
    \mathbf{E}(z;x,t) = \begin{pmatrix} 1&  0 \\ -\overline{(\mathbf{E}_1(x,t))_{12}}  & 1 \end{pmatrix} + \frac{\mathbf{E}_1(x,t)}{z} + \bigO{z^{-2}} ,\\
    2i(\mathbf{E}_1)_{12} = \frac{1}{\sqrt{2t}} 
      \left[ 
        \beta_{12} \mathbf{P}^\textrm{sol}_{11}(\lambda_0;x,t)^2 
       + \beta_{21} \mathbf{P}^\textrm{sol}_{12}(\lambda_0;x,t)^2	
      \right] 
      + \bigO{t^{-1}} .
  \end{gathered}
\end{equation}

\begin{remark}
	Though we simply ignored the $\overline{\partial}$-derivatives of $\mathbf{N}^{(2)}$ to arrive at the model problem $\mathbf{P}$, it turns out that both the leading order solitonic component of the solution and the first dispersive correction terms are encoded in the solution of $\mathbf{P}$. In particular, the dispersive correction, is up to a trivial multiplicative constant, is given exactly by $2i(\mathbf{E}_{1})_{12}$.
\end{remark}

\subsection{Step 4: Estimating the $\overline{\partial}$ contribution as $t \to \infty$.} 
We now return to our non-analytic, piecewise continuous unknown $\mathbf{N}^{(2) }$ defined by \eqref{N2}. The solution $\mathbf{P}(z;x,t)$ of our model problem, Problem~\ref{model}, was defined such that it has exactly the same jump discontinuities as $\mathbf{N}^{(2)}$ across each component of $\Gamma^{(2)}$ (cf \eqref{Lambda.2}. We use $\mathbf{P}$ to make one further transformation which has the effect of removing the jump discontinuities from the problem. Let
\begin{equation}\label{N3} 
	\mathbf{N}^{(3)}(z;x,t) = \mathbf{N}^{(2)} (z;x,t) \mathbf{P}(z;x,t)^{-1}.
\end{equation}
As $\mathbf{N}^{(2)}$ is continuous up to the boundary of each connected compent of $\C \setminus \Gamma^{(2)}$ and has the same jump conditions \eqref{P.jump}-\eqref{P.soliton} as $\mathbf{P}$ across each component, it follows that 
\[
\begin{multlined}[b][.65\displaywidth]	
	\mathbf{N}^{(3)}_+ (\lambda;x,t) 
	= \mathbf{N}^{(2)}_+ (z;x,t) \mathbf{P}_+(z;x,t)^{-1} \\
	= \mathbf{N}^{(2)}_- (z;x,t) \mathbf{J}_{x,t}^{(P)}(\lambda) 
	  \left( \mathbf{P}_-(z;x,t) \mathbf{J}_{x,t}^{(P)}(\lambda) \right)^{-1} 
	= \mathbf{N}^{(3)}_- (\lambda;x,t)  
\end{multlined} \qquad  \lambda \in \Gamma^{(2)}.
\]
is continuous across each of the jump boundaries $\Gamma^{(2)}$. However, it is still non-analytic due to the extension \eqref{extensions} used in \eqref{N2} to define $\mathbf{N}^{(2)}$. It can be shown that the new unknown $\mathbf{N}^{(3)}$ satisfies the following problem.
%
%
\begin{problem}\label{dbar} Let $x \in \R$, $t>0$ be parameters. Find a $1\times 2$ vector function $\mathbf{N}^{(3)} = \mathbf{N}^{(3)}(u,v) = \mathbf{N}^{(3)}(u,v;x,t)$, $(u,v) \in \R^2$ with the following properties: 
\begin{itemize}
	\item[(i)] (Continuity) $\mathbf{N}^{(3)}$ is a continuous function of $(u,v) \in \R^2$.
	\item[(ii)] (Normalization) $\mathbf{N}^{(3)}(u,v) \to \I$ as $(u,v) \to \infty$.
	\item[(iii)] (Nonanalyticity) $\mathbf{N}^{(3)}$ is a (weak) solution of the partial differential equation $\dbar \mathbf{N}^{(3)}(u,v) = \mathbf{N}^{(3)}(u,v) \mathbf{W}(u,v)$ where $\mathbf{W}(u,v) = \mathbf{W}(u,v;x,t)$ is defined as follows.
	\begin{equation}
	\begin{gathered}
		\mathbf{W}(u,v;x,t) = \mathbf{P}(u+iv;x,t) \dbar \mathbf{R}(u,v;x,t) \mathbf{P}(u+iv;x,t)^{-1} \\
		\dbar \mathbf{R}(u,v;x,t) = \begin{cases}
		\twomat{0}{0}{-\dbar R_1(u,v) e^{2it \theta(u+iv,x/t) } }{0} & u+iv \in \Omega_1 \\
		\\
		\twomat{0}{-\dbar R_3(u,v) e^{-2it \theta(u+iv, x/t) } }{0}{0} & u+iv \in \Omega_3 \\
		\\
		\twomat{0}{0}{\dbar R_4(u,v) e^{2it \theta(u+iv, x/t) } }{0} & u+iv \in \Omega_4 \\
		\\
		\Twomat{0}{\dbar R_6(u,v) e^{-2it \theta(u+iv, x/t) } }{0}{0} & u+iv \in \Omega_6 \\
		\\
		\mathbf{0} & u+iv \in \Omega_2 \cup \Omega_5
		\end{cases}
	\end{gathered}
	\end{equation}
\end{itemize}
\end{problem}

The solution of this $\dbar$-problem can be represented as a Fredholm integral equation over $\R^2$  
\begin{gather}\label{n3.int.eq}
	(1 - K_W)  \mathbf{N}^{(3)} = (1,\ 0) \\
\nonumber	
(K_W f)(u,v;x,t) = 
	 \frac{1}{\pi} \int_{\R^2} \frac{f(u', v' ; x,t) \mathbf{W}(u',v';x,t) }{ (u+iv) - (u'+iv')} du' dv'.
\end{gather}

The following proposition completes the construction of the solution of our inverse problem. 
\begin{proposition}[\cite{JLPS18b}, Proposition 3.18]\label{prop:N3}
Suppose that $\rho \in H^{2,2}(\R)$ with $c \coloneqq \inf_{\lam \in \R} \left(1- \lam |\rho( \lam)|^2 \right) >0$ strictly.
Then, for sufficiently large time $t>0$, there exists a unique solution $\mathbf{N}^{(3)}(u,v;x,t)$ for $\dbar$-Problem~\ref{dbar} with the property that 
\begin{equation}
\label{N3.exp}
\mathbf{N}^{(3)}(u,v;x,t) = (1\ 0) + \frac{1}{u+iv}  \mathbf{N}^{(3)}_1(x,t) + o_{\xi,t} \left( \frac{1}{|u+iv|} \right),
\quad (u,v) \to \infty
\end{equation}
where
\begin{equation}
\label{N31.est}
\left| \mathbf{N}_1^{(3)}(x,t) \right| \lesssim t^{-3/4}
\end{equation}
 where the implied constant in \eqref{N31.est} is independent of $\xi$ and $t$ and uniform for 
$\rho$ in a bounded subset of $H^{2,2}(\R)$ with $\inf_{\lam \in \R} (1-\lam|\rho(\lam)|^2) \geq c>0$ for a fixed $c>0$.
\end{proposition}

The idea of the proof is as follows. As $\mathbf{P}$ is bounded and boundedly invertible, 
we may use the fact\footnote{See, for example, the book of Astala, Iwaniec, and Martin \cite[Section 4.3]{AIM09}} the the solid Cauchy transform
\[
	(P f )(u,v) = \frac{1}{\pi}\int_{\R^2} \frac{1}{u-u' + i (v-v') } f(u',v') du'dv'
\]
is bounded from $L^p(\R^2) \cap L^q(\R^2)$ to $L^\infty(\R^2)$ for any $p,q$ with $1 < p < 2 < q < \infty$ and Lemma \ref{lem:R.dbar} to bound the operator $K_W$. Indeed, it's shown in \cite{JLPS18b} that $\| K_W \|_{L^\infty(\R^2) \to L^\infty(\R^2)} \lesssim t^{-1/4}$. This allows us to establish the existence of a unique solution to \eqref{n3.int.eq} and to compute the asymptotic expansion in Proposition~\ref{prop:N3} by asymptotically expanding the resolvent operator 
 $(\I-K_W)^{-1}$ in a Neumann series for large times.

\subsection*{Large-time asymptotic behavior of $q(x,t)$}
Inverting the sequence of transformation from $\mathbf{N}$ to $\mathbf{N}^{(3)}$, the solution of Riemann-Hilbert Problem~\ref{RHP longtime} is given for large $z$ by
\[
	\mathbf{N}(z) = \mathbf{N}^{(3)}(z) \mathbf{E}(z) \mathbf{P}^\mathrm{sol}(z)\mathbf{R}(z)^{-1} \delta(z,\lambda_0)^{\sigma_3}.
\]
It follows from \eqref{delta.expand} and Lemma~\ref{lem:R.dbar} that 
\[
\delta(z)^{\sigma_3} = \I + (\delta_1/z) \sigma_3 + \bigO{z^{-2}},
\qquad
\mathbf{R}(z) = \I + \bigO{e^{-ct} }.
\]
Inserting these, Proposition~\ref{prop:Nsol} and the expansion \eqref{model.error} and \eqref{N3.exp} into the recovery formula \eqref{DNLS2:q.recon.bis}  for $q(x,t)$ we have
\[
	q(x,t) = q_\mathrm{sol}(x,t; \mathcal{D}_\xi) + 2i (\mathbf{E}_1)_{12} + \bigO{t^{-3/4}}.
\]
Here,  $2i(\mathbf{E}_1)_{12}$ gives, up to the gauge transformation $u = \mathcal{G}^{-1}(q)$, the leading order behavior of the dispersive correction in Theorem~\ref{thm:soliton-resolution}.

To recover the  large-time asymptotic behavior of the solution $u(x,t)$ of \eqref{DNLS} one needs to compute an asymptotic expansion of the gauge factor $\exp \left( i \int_{x}^\infty  |q(y,t)|^2 dy \right)$ relating the two solutions. A straightforward computation shows that this gauge factor can be computed directly from the solution of the Riemann-Hilbert problem in the form $\exp \left( i \int_{x}^\infty  |q(y,t)|^2 dy \right) = \mathbf{N}_{11}(0;x,t)^{-2}$. By using the asymptotic expansion of $\mathbf{N}$ we have constructed above, an asymptotic expansion for the gauge factor can be computed\footnote{There is an extra complication in this computation when $\lambda_0 \to 0$---the gauge factor then depends on the local model $\mathbf{P}^\textrm{pc}(\zeta(0);x,t)$ in a non-trivial way--which is why we exclude the region $|\lambda_0| < Mt^{-1/8}$ from Theorem~\ref{thm:soliton-resolution}}, and one recovers the asymptotic expansions for $u(x,t)$ described by Theorem~\ref{thm:soliton-resolution}.

\section*{Acknowledgements}

This work was supported by a grant from the Simons Foundation/SFARI (359431, PAP).
CS  is supported in part by  Discovery Grant 2018-04536 from the Natural Sciences and Engineering Research Council of Canada.

%
%

\end{document}